\documentclass{amsart}
\usepackage{amsmath, amsthm, amssymb}
\usepackage{enumerate}
\usepackage{geometry}
\usepackage{xcolor}
\usepackage[colorlinks=true, citecolor=red, linkcolor=blue, urlcolor=black]{hyperref}

\newtheorem{theorem}{Theorem}[section]
\newtheorem{proposition}{Proposition}[section]
\newtheorem{lemma}{Lemma}[section]
\newtheorem{remark}{Remark}[section]
\newtheorem{definition}{Definition}[section]

\numberwithin{equation}{section}

\newcommand{\R}{\mathbb{R}}

\allowdisplaybreaks

\title[Weakly competing bubbles for a critical elliptic system]
{Weakly competing bubbles for a critical elliptic system in higher dimensions}

\author[I. Ceresa Dussel]{Ignacio Ceresa Dussel}
\address{Ignacio Ceresa Dussel
\newline\indent Departamento de Matem\'aticas, Universidad Aut\'onoma de Madrid, Madrid, Spain,
\newline\indent Instituto de Investigaciones Matem\'aticas ``Luis A. Santal\'o'' (IMAS), UBA--CONICET, Buenos Aires, Argentina}
\email{juan.ceresa@uam.es}

\author{Maria Medina}
\address{Maria Medina
\newline\indent Departamento de Matem\'aticas, Universidad Aut\'onoma de Madrid, Madrid, Spain
}
\email{maria.medina@uam.es}

\author[A. Pistoia]{Angela Pistoia}
\address{Angela Pistoia
\newline\indent Dipartimento di Scienze di Base e Applicate per l'Ingegneria, Sapienza Universit\`a di Roma, Rome, Italy}
\email{angela.pistoia@uniroma1.it}


\begin{document}

\begin{abstract}
We consider the critical competitive system
\begin{equation*}
\begin{cases}
-\Delta u=u^{\frac{n+2}{n-2}}
+\beta u^{\frac{2}{n-2}}v^{\frac{n}{n-2}}
&\text{in }\R^n,\\[2pt]
-\Delta v=v^{\frac{n+2}{n-2}}
+\beta v^{\frac{2}{n-2}}u^{\frac{n}{n-2}}
&\text{in }\R^n,
\end{cases}
\end{equation*}
where
$p=\frac{n+2}{n-2}$, $q=\frac{2}{n-2}$ and $\beta<0$.
In dimensions three and four, solutions combining a bubble in one
component with a ring of bubbles in the other have been constructed
by Lyapunov--Schmidt reduction. For $n\ge5$ this approach breaks down:
the coupling is sublinear, the energy is only of class $C^1$, and small
corrections cannot be kept away from the region where the nonlinearity
is singular. We overcome this obstruction with a variational reduction
in the spirit of Bahri and Rey, in which the infinite-dimensional
correction is obtained by constrained minimization and is only required
to depend continuously on the parameters. The resulting
finite-dimensional problem is then solved by a degree argument.
In this way, for every $n\ge5$ and every integer $k\ge2$, we construct
nonnegative finite-energy solutions in the weak-competition regime
$\beta\to0^-$: the first component is close to a standard bubble, while
the second concentrates at the $k$ vertices of a regular polygon. The
concentration scale is $d_*\varepsilon(1+o(1))$, with $d_*$ explicit and
$\varepsilon^{(n-4)/2}=|\beta|\,|\ln\varepsilon|$, and we give precise
expansions of the amplitudes and of the remainders.
\end{abstract}

\maketitle

\section{Introduction}

In this paper we study the critical elliptic system
\begin{equation}\label{eq:system}
\begin{cases}
-\Delta u=u^{\frac{n+2}{n-2}}
+\beta u^{\frac{2}{n-2}}v^{\frac{n}{n-2}}
&\text{in }\R^n,\\[2pt]
-\Delta v=v^{\frac{n+2}{n-2}}
+\beta v^{\frac{2}{n-2}}u^{\frac{n}{n-2}}
&\text{in }\R^n,\\
u,v\ge0,\qquad u,v\in\dot H^1(\R^n),
\end{cases}
\end{equation}
where $n\ge3$ and $\beta\in\R$ is a coupling parameter.

Systems of this kind appear in the study of standing waves for coupled
nonlinear Schr\"odinger (Gross--Pitaevskii) equations, which describe,
for instance, binary mixtures of Bose--Einstein condensates
\cite{EsryGreeneBurkeBohn1997,Timmermans1998} and the propagation of
mutually incoherent light beams in nonlinear optics
\cite{AkhmedievAnkiewicz1999}. Each component represents the density
amplitude of one species, while $\beta$ measures the interaction between
the two species. When $\beta>0$ the interaction is attractive, or
\emph{cooperative}, and the two species tend to coexist; when $\beta<0$
it is repulsive, or \emph{competitive}, and they tend to separate in
space \cite{Timmermans1998}. System \eqref{eq:system} can be regarded
as a higher-dimensional, conformally invariant counterpart of the cubic
Gross--Pitaevskii model, in which all the nonlinearities have critical
Sobolev growth.

For $\beta=0$ the two equations decouple, and each component solves the
Yamabe equation $-\Delta w=w^{(n+2)/(n-2)}$ in $\R^n$, whose positive
finite-energy solutions are the standard bubbles
\cite{Aubin1976,Talenti1976,CaffarelliGidasSpruck1989}
$$
U_{\delta,\xi}(x)
=C_n\left(\frac{\delta}{\delta^2+|x-\xi|^2}\right)^{\frac{n-2}{2}},
\qquad
C_n=[n(n-2)]^{\frac{n-2}{4}},
\qquad \delta>0,\quad\xi\in\R^n.
$$
Thus, for every $\beta$, system \eqref{eq:system} has the
\emph{semi-trivial} solutions $(U_{\delta,\xi},0)$ and
$(0,U_{\delta,\xi})$, and the interesting question is the existence of
\emph{fully nontrivial} solutions, namely solutions with both components
different from zero. The simplest ones are the \emph{synchronized}
solutions, in which both components are multiples of the same bubble:
for every $\beta>-1$,
$$
u=v=(1+\beta)^{-\frac{n-2}{4}}U_{\delta,\xi}
$$
solves \eqref{eq:system}. In the cooperative regime nothing else can
happen: Guo and Liu \cite{GuoLiu2008} proved that, for $\beta>0$, every
positive solution of \eqref{eq:system} is synchronized. Non-synchronized
solutions must therefore be sought in the competitive regime $\beta<0$,
which is the setting of this paper.

When $\beta<0$, positive solutions of the system are closely related to
sign-changing solutions of the single equation
\begin{equation}\label{eq:nodal}
-\Delta w=|w|^{\frac{4}{n-2}}w\quad\text{in }\R^n .
\end{equation}
Roughly speaking, as $\beta\to-\infty$ the two components of suitable
solutions $(u_\beta,v_\beta)$ separate from each other, and
$u_\beta\sim w^+$, $v_\beta\sim w^-$, where $w$ is a sign-changing
solution of \eqref{eq:nodal}. This phase separation phenomenon was
first discovered in the subcritical regime by Terracini and her
collaborators
\cite{ContiTerraciniVerzini2002,ContiTerraciniVerzini2003,
ContiTerraciniVerzini2005,NorisTavaresTerraciniVerzini2010,
TavaresTerracini2012,SoaveZilio2015,SoaveTavaresTerraciniZilio2016}
and by Caffarelli and Lin \cite{CaffarelliLin2008}. In the critical
setting, it has been established in bounded domains by Chen
and Zou \cite{ChenZou2012,ChenZou2015} (see also
\cite{TavaresYouZou2022}), and for the Yamabe equation on the sphere and
on closed manifolds in \cite{ClappSaldanaSzulkin2021,ClappPistoiaTavares}. It is then natural
to ask the converse question: {\em given a sign-changing solution $w$ of
\eqref{eq:nodal} with $\ell$ nodal domains $\Omega_1,\dots,\Omega_\ell$,
does the competitive system with $\ell$ components have a positive
solution whose $i$-th component converges to the restriction of $|w|$
to $\Omega_i$ as $\beta\to-\infty$?}

In this generality the question seems out of reach. A more feasible
approach is to mimic, at the level of the system, the methods that have
been used to build sign-changing solutions of the single equation. On
the variational side, Ding \cite{Ding1986} found infinitely many
sign-changing solutions of \eqref{eq:nodal} which are invariant under
the conformal action of the group $O(m)\times O(n+1-m)$ induced by the
stereographic projection. Positive solutions of critical systems with
several components and the same symmetries were later obtained by Clapp
and Pistoia \cite{ClappPistoia2018}, Clapp and Szulkin
\cite{ClappSzulkin2019}, and Clapp, Salda\~na and Szulkin
\cite{ClappSaldanaSzulkin2021}. 

On the perturbative side, several families of sign-changing solutions
of \eqref{eq:nodal} have been constructed in
\cite{delPinoMussoPacardPistoia2011,delpino2013torus,
medina2019desingularization,MedinaMusso2021,LiSun2026}; low-energy
nodal solutions have also been obtained by ODE techniques in
\cite{FernandezPetean2020}. Our starting point is the \emph{crown}
solution of del Pino, Musso, Pacard and Pistoia
\cite{delPinoMussoPacardPistoia2011}: for every $k$ large enough,
\eqref{eq:nodal} has a solution of the form
$$
w(y)=U(y)-\sum_{\ell=1}^k\delta^{-\frac{n-2}{2}}
U\Big(\frac{y-\xi_\ell}{\delta}\Big)+\phi(y),
$$
where $U:=U_{1,0}$, $\delta\to0$ and $\phi\to0$ as $k\to\infty$, and the
points $\xi_\ell$ are close to the vertices of a regular $k$-polygon
inscribed in the unit circle of the $(y_1,y_2)$-plane. In other words,
$w$ looks like a positive bubble surrounded by many negative bubbles.
It is then natural to look for solutions of \eqref{eq:system} in which
the positive and the negative parts of the crown are carried by
different components, that is,
$$
u\approx U\qquad\text{and}\qquad
v\approx\sum_{\ell=1}^kU_{\delta,\xi_\ell}.
$$
Solutions of this type have been found by Fern\'andez, Medina and
Pistoia \cite{FernandezMedinaPistoia2026} in dimension three, in the
strong competition regime $\beta\to-\infty$, with
$\delta_\beta^{1/2}|\ln\delta_\beta|\sim|\beta|^{-1}$, and by Chen,
Medina and Pistoia \cite{ChenMedinaPistoia2023} in dimension four, in
the weak competition regime $\beta\to0^-$, with
$|\beta|\,|\ln\delta_\beta|\sim1$. Unlike in
\cite{delPinoMussoPacardPistoia2011}, in both cases the number $k\ge2$
of bubbles is fixed, and the small parameter is the coupling strength.
These results show that the dimension determines the coupling regime in
which this configuration appears. In this paper we deal with all
dimensions $n\ge5$: the construction again takes place in the weak
competition regime $\beta\to0^-$.

To state our result, let $X_k\subset\dot H^1(\R^n)$ be the space of
functions which are even in each of the variables $x_2,\dots,x_n$,
invariant under the rotation of angle $2\pi/k$ in the $(x_1,x_2)$-plane,
and invariant under the Kelvin transform (see Section~\ref{sec:setting}).

\begin{theorem}\label{thm:main-m2}
Let $n\ge5$ and let $k\ge2$ be an integer.
There exists $\beta_0<0$ such that, for every
$\beta\in[\beta_0,0)$, system \eqref{eq:system}
admits a solution
$$
(u_\beta,v_\beta)\in X_k\times X_k
\subset\dot H^1(\R^n)\times\dot H^1(\R^n)
$$
of the form
$$
u_\beta=a_\beta U+\phi_\beta,
\qquad
v_\beta=b_\beta\sum_{j=1}^k
U_{\delta_\beta,\xi_j}+\psi_\beta,
$$
where
$$
\xi_{j}=\rho_\beta
\left(\cos\frac{2\pi(j-1)}k,\sin\frac{2\pi(j-1)}k,0,\ldots,0\right),
\qquad \rho_\beta^2+\delta_\beta^2=1.
$$
Let $\varepsilon\in (0,1)$ be given by
$$
\varepsilon^{\frac{n-4}{2}}=|\beta|\,|\ln\varepsilon|.
$$
As $\beta\to0^-$, the parameters satisfy $a_\beta, b_\beta\to 1$ and
$$
\delta_\beta=d_*\varepsilon+o(\varepsilon),\qquad
\|\phi_\beta\|_{\dot H^1}+\|\psi_\beta\|_{\dot H^1}
=O\!\left(\frac{\varepsilon^{(n-2)/2}}{|\ln\varepsilon|}\right),
$$
where $d_*>0$ is an explicit constant, depending only on $n$ and $k$.
\end{theorem}

Although in every dimension the problem is eventually reduced to a
finite-dimensional one, the way this reduction is carried out depends
strongly on $n$. Set $q:=2/(n-2)$. The coupling term in the first
equation is $\beta u^qv^{q+1}$, which is superlinear in $u$ when $n=3$,
linear when $n=4$ and \emph{sublinear} when $n\ge5$; the same holds for
the second equation, with $u$ and $v$ exchanged. In dimensions three and
four the classical Lyapunov--Schmidt reduction can be applied
\cite{FernandezMedinaPistoia2026,ChenMedinaPistoia2023}. For $n\ge5$,
instead, $0<q<1$: the map $s\mapsto s_+^q$ is not locally Lipschitz at
the origin, and the coupling term in the energy is of class $C^1$ but
not $C^2$. Moreover, a correction which is small in $\dot H^1(\R^n)$
need not be small pointwise compared with the bubbles, so there is no
way to keep the perturbed components away from the region where the
coupling is singular. For these reasons, the standard
Lyapunov--Schmidt reduction, based on a contraction mapping or on a
smooth implicit function theorem, cannot be applied directly.

The key idea of this paper is to replace it with a variational
reduction in the spirit of Bahri, Rey, and Bahri--Li--Rey
\cite{Bahri1989,Rey1990,BahriLiRey1995}. Rather than solving the
infinite-dimensional equations by a fixed point argument, we obtain the
correction by minimizing the energy in a small, suitably constrained
neighborhood. The advantage is that we do not need the correction to
depend differentiably on the finite-dimensional parameters: continuity
is enough for the final degree argument. In other words, the
variational nature of the reduction is not just a different way of
presenting the usual scheme, but precisely what allows us to deal with
the non-Lipschitz coupling.

To carry out this program, we combine rotational and reflection
symmetries with the Kelvin invariance. Once the amplitude directions  
and the remaining approximate kernel direction are removed, the
uncoupled quadratic form is uniformly coercive. For the first component
this follows from the nondegeneracy of the standard bubble
\cite{Rey1990,BianchiEgnell1991}. For the concentrating component we use
a blow-up argument together with the spectral gap of a single bubble;
here the Kelvin invariance rules out the last translation mode. The coupling energy is not
jointly convex, but its contribution to the monotonicity estimate is
bounded from below by $-C|\beta|$ times a quadratic term. Hence, for
$|\beta|$ small, the restricted energy is strictly convex in a small
ball, and the correction exists, is unique, and is obtained as the
minimum point of the energy in that ball.

Finally, the correction depends continuously on the parameters.
Testing the reduced equation against the three remaining directions gives a
continuous system in the two amplitudes $a$ and $b$ and the concentration scale $\delta$
which, after rescaling, converges uniformly to an explicit vector field
with a nondegenerate zero. A Brouwer degree argument then provides a
zero of the reduced system.

\begin{remark}\rm
The proof also gives the first-order corrections of the amplitudes
$a_\beta$ and $b_\beta$, which are of order $\varepsilon^{n-2}$; see
\eqref{eq:amplitudes}. In particular, $a_\beta$ is slightly larger and
$b_\beta$ slightly smaller than $1$. Hence the first component converges
strongly in $\dot H^1(\R^n)$ to $U$, whereas the second one is close, in
the energy norm, to $k$ bubbles concentrating at the vertices of a
regular polygon inscribed in the unit circle.

The concentration rate comes from the balance of two interactions. If
$V_\delta:=\sum_{j=1}^kU_{\delta,\xi_j}$, the interaction among the
bubbles of $V_\delta$ is of order $\delta^{n-2}$, while the interaction
between $U$ and $V_\delta$ is of order
$$
|\beta|\int_{\R^n}U^{\frac n{n-2}}V_\delta^{\frac n{n-2}}\,dx
\asymp |\beta|\,\delta^{n/2}|\ln\delta|.
$$
The logarithm is due to the borderline decay of $U^{n/(n-2)}$ after
rescaling around a concentration point: away from $\xi$, the function
$U_{\delta,\xi}^{n/(n-2)}$ decays exactly like $|x-\xi|^{-n}$.
Balancing the two terms gives $\delta^{(n-4)/2}\asymp|\beta|\,|\ln\delta|$,
which explains the definition of $\varepsilon$; for $n\ge5$, this forces
$\delta\to0$ as $\beta\to0^-$. In the reduced equation for the scale,
the two leading terms have opposite signs, and their balance determines
the positive constant $d_*$ given by the identity \eqref{dvalue}.
\end{remark}

\begin{remark}\rm
In dimensions three and four, the constructions in
\cite{FernandezMedinaPistoia2026,ChenMedinaPistoia2023} also apply to
systems with one component $u$ and $m$ components $v_1,\dots,v_m$,
which compete with $u$ and interact with each other through an
additional coupling parameter. In those solutions each $v_i$ is a
rotation of $v_1$, so that the concentration points lie at the
vertices of rotated copies of the same polygon. It would be interesting
to know whether such configurations also exist for $n\ge5$, and whether
the full symmetry among the components $v_i$ can be broken.
\end{remark}

\begin{remark}\rm
Comparing dimension three with higher dimensions, one may also ask
whether, for $n\ge4$, there are other configurations of concentration
points which give rise to solutions in the strong competition regime
$\beta\to-\infty$.
\end{remark}

The paper is organized as follows. In Section~\ref{sec:setting} we
introduce the ansatz, the symmetric function spaces and the projected
equations. Section~\ref{sec:infinite-dim} contains the error and
coercivity estimates and the construction of the correction by
minimization. In Section~\ref{sec:finite-dim} we derive the reduced
equations and compute the interaction terms. Finally, in the last
section we solve the reduced system by a degree argument and complete
the proof of Theorem~\ref{thm:main-m2}.

\section{The problem: ansatz and setting}\label{sec:setting}

Let $n\ge5$, $k\ge2$ be a fixed integer, and $\beta<0$, and consider the problem
\begin{equation*}
    \left\{
    \begin{aligned}
        &-\Delta u = u^{\frac{n+2}{n-2}} + \beta\, u^{\frac{2}{n-2}} v^{\frac{n}{n-2}}
        && \text{in } \mathbb{R}^n, \\
        &-\Delta v = v^{\frac{n+2}{n-2}} + \beta\, v^{\frac{2}{n-2}} u^{\frac{n}{n-2}}
        && \text{in } \mathbb{R}^n.
    \end{aligned}\right.
\end{equation*}
Since we are seeking nonnegative solutions, it is convenient to consider the following
system involving the positive and negative parts of a function, denoted by
$s_+ := \max\{s,0\}$ and $s_- := \min\{s,0\}$:
\begin{equation}\label{eq:positive-system}
    \left\{
    \begin{aligned}
        &-\Delta u = u_{+}^{\frac{n+2}{n-2}} + \beta\, u_{+}^{\frac{2}{n-2}} v_{+}^{\frac{n}{n-2}}
        && \text{in } \mathbb{R}^n, \\
        &-\Delta v = v_{+}^{\frac{n+2}{n-2}} + \beta\, v_{+}^{\frac{2}{n-2}} u_{+}^{\frac{n}{n-2}}
        && \text{in } \mathbb{R}^n.
    \end{aligned}\right.
\end{equation}
If $(u,v)$ is a solution to \eqref{eq:positive-system}, then by testing with $u_-$ and $v_-$,
one readily obtains $u_- \equiv v_- \equiv 0$ in $\R^n$. Consequently, $(u,v)$ is
nonnegative and therefore also a solution to the original system. Hence, we will focus
on solving \eqref{eq:positive-system}.

We recall that the critical Sobolev exponent is $2^{*} := \frac{2n}{n-2}$, and set
$p := 2^{*} - 1 = \frac{n+2}{n-2}$.
Let us consider the  single critical equation
$$
    -\Delta U(x) = U^p(x) \quad \text{in } \R^n.
$$
It is well known that all nontrivial nonnegative solutions in $\dot H^1(\R^n)$ are the {\em bubbles}
\begin{equation}\label{eq:standard-bubble2}
    U_{\delta,\,\xi}(x)
    := \delta^{-\frac{n-2}{2}}\,
    U\!\left( \frac{x - \xi}{\delta} \right),\quad \delta>0,\ \xi\in\mathbb R^n
\end{equation}
where
\begin{equation}\label{eq:standard-bubble}
    U(x) := \frac{C_n}{\left( 1 + |x|^2 \right)^{\frac{n-2}{2}}},
    \qquad \text{where } C_n = \left( n(n-2) \right)^{\frac{n-2}{4}}.
\end{equation}
We look for a solution to \eqref{eq:positive-system} of the form
$$
    u = aU + \phi, \qquad v = bV + \psi,
$$
where
$$
   V(x) := \sum_{j=1}^k U_j(x),\quad   U_j(x) := U_{\delta,\,\xi_j}(x).
$$
The concentration points $\xi_j$
in $\R^n$  are arranged at the vertices of a regular polygon, that is,
\begin{equation}\label{eq:vertices}
    \xi_j := \rho\left( \cos\frac{2\pi(j-1)}{k},\,\sin\frac{2\pi(j-1)}{k},\,
    0,\ldots,0\right),
    \qquad \rho^2 + \delta^2 = 1,
    \qquad j = 1,\ldots,k.
\end{equation}
Moreover, we choose the parameters in the form
\begin{equation}\label{eq:rescaling}
a=1+t\varepsilon^{n-2},\qquad
b=1+s\varepsilon^{n-2},\qquad
\delta=d\varepsilon,\qquad d>0,
\end{equation}
where $\varepsilon\in(0,1)$ is uniquely determined, for $\beta<0$, by
\begin{equation}\label{eq:eps-beta}
\varepsilon^{\frac{n-4}{2}}=|\beta|\,|\ln\varepsilon|.
\end{equation}
Indeed, the map $\sigma\mapsto\sigma^{-\frac{n-4}{2}}|\ln\sigma|$
is strictly decreasing from $+\infty$ to $0$ on $(0,1)$, so that
$\varepsilon\to0$ as $\beta\to0^-$. 
The remainder terms 
$
    \phi \in  \mathcal{H}_{1}$ and $ \psi \in  \mathcal{H}_{2}.
$
The spaces $ \mathcal{H}_{1}$ and $ \mathcal{H}_{2}$ are defined in the next section.

We will work in the homogeneous Sobolev space $\dot{H}^1(\R^n)$, defined as the
completion of $C_c^\infty(\R^n)$ with respect to the norm
$$
\|u\| := \left( \int_{\R^n} |\nabla u|^2 \, dx \right)^{1/2}.
$$
We denote by $\langle \cdot,\cdot\rangle_{H}$ the associated inner product,
$$
\langle f, g\rangle_{H} := \int_{\R^n} \nabla f \cdot \nabla g \, dx.
$$
We reserve $\langle\cdot,\cdot\rangle_{H}$ for this inner product, $\langle\cdot,\cdot\rangle$ always denotes a duality pairing between
$\dot H^1(\R^n)$ and its dual, as in \eqref{eq:residuals} below.

\begin{definition}We denote by $X_k$ the space of functions $h\in \dot{H}^{1}(\R^n)$ satisfying:
\begin{enumerate}
    \item $h$ is even in  {each of the variables $x_2,\dots,x_n$}.
    \item $h$ is invariant under rotations of angle $\frac{2\pi}{k}$ in the
    $(x_1, x_2)$-plane.
    \item $h$ is invariant under the Kelvin transform:
    $$
    h(x) = \frac{1}{|x|^{n-2}} h\!\left( \frac{x}{|x|^2} \right).
    $$
\end{enumerate}
\end{definition}
It is worthwhile to observe that both $U$ and $V$ belong to $X_k$.

Let us  consider the linearization of the Yamabe equation around $U$, i.e.,
\begin{equation}\label{eq:linearization}
    -\Delta\phi(x) = p\,U(x)^{p-1}\phi(x) \quad \text{in } \R^n.
\end{equation}
The set of $\dot{H}^1$ solutions to \eqref{eq:linearization} is
$\mathrm{span}\{Z^i : i=0,1,\ldots,n\}$, where $Z^0$ and $Z^i$ are the dilation and
translation modes of $U$,
\begin{equation}\label{Zi}\begin{aligned}
&Z^0:=\delta\,\partial_\delta U_{\delta,0}\big|_{\delta=1}
=\frac{n-2}{2}\,C_n\,\frac{|x|^2-1}{(1+|x|^2)^{n/2}},
\\
&Z^i:=\partial_{\xi_i}U_{1,\xi}\big|_{\xi=0}
=(n-2)\,C_n\,\frac{x_i}{(1+|x|^2)^{n/2}},\quad i=1,\dots,n .\end{aligned}
\end{equation}
Set
\begin{equation}\label{eq:Zdelta}
    Z_\delta := \delta\,\frac{d}{d\delta}V=\sum_{j=1}^k Z_j,
    \qquad
    Z_j:=\delta\,\frac{d}{d\delta}\big[U_{\delta,\xi_j}\big]
    =\delta\, \partial_\delta U_{\delta,\xi_j}
    - \frac{\delta^2}{\rho}\, \zeta_j \cdot \nabla_\xi U_{\delta,\xi_j},
    \qquad \zeta_j := \frac{\xi_j}{\rho}.
\end{equation}
Let us introduce the finite-dimensional  spaces 
$$
\mathcal{F}_1 := \operatorname{span}\{U\},
\qquad
\mathcal{F}_2 := \operatorname{span}\{V, Z_\delta\},
$$
and their orthogonal spaces
$$
\mathcal{H}_{1} := \{\phi \in X_k : \langle \phi, U\rangle_{H} = 0\},
\qquad
\mathcal{H}_{2} :=
 \left\{ \psi \in X_k : \langle \psi, V\rangle_{H} = 0 \ \text{and} \ \langle \psi, Z_\delta\rangle_{H} = 0 \right\}.
$$
It can be checked that both $\mathcal{F}_1, \mathcal{F}_2$ are contained in $X_k$. Thus, by construction, we have the orthogonal decompositions
$$
X_k = \mathcal{F}_1 \oplus  \mathcal{H}_{1},
\qquad
X_k = \mathcal{F}_2 \oplus  \mathcal{H}_{2}.
$$
Let $P_1: X_k \to \mathcal{F}_1$ and $P_2: X_k \to \mathcal{F}_2$ denote the
associated orthogonal projections, and let
$$
P_1^\perp := \mathrm{Id} - P_1 : X_k \longrightarrow  \mathcal{H}_{1},
\qquad
P_2^\perp := \mathrm{Id} - P_2 : X_k \longrightarrow  \mathcal{H}_{2}
$$
be their complements. Notice that $\mathcal{F}_2$ and $\mathcal{H}_{2}$ depend on $\delta$.

By the Sobolev inequality, every element of $\dot{H}^1(\R^n)$ can be identified with a function
in $L^{2^*}(\R^n)$ whose distributional gradient belongs to $L^2(\R^n)$. We
define, via the Riesz representation theorem, the operator
$(-\Delta)^{-1} : L^{(2^*)'}(\R^n) \to \dot{H}^1(\R^n)$ as
$$
(-\Delta)^{-1}(f) = u \iff (-\Delta)u = f \quad \text{in } \R^n.
$$
We use the same notation for the Riesz inverse on the full dual
$(\dot H^1(\R^n))^*$. Throughout we write
$$
    q := \frac{p-1}{2} = \frac{2}{n-2},
    \qquad\text{so that}\qquad
    q+1 = \frac{p+1}{2} = \frac{n}{n-2}.
$$
 {Since for $u,v\in\dot H^1(\R^n)$ the nonlinear terms $u_+^p$, $u_+^qv_+^{q+1}$ (and
the corresponding ones in the second equation) only belong to $L^{(2^*)'}(\R^n)$, it is
convenient to formulate \eqref{eq:positive-system} weakly.}
For $u,v \in \dot H^1(\R^n)$ we therefore define the \emph{residual functionals}
by duality,
\begin{equation}\label{eq:residuals}
\begin{aligned}
    \langle R_1(u,v), h\rangle
    &:= \int_{\R^n}\Big( \nabla u\cdot\nabla h - u_+^{p}\, h
        - \beta\, u_+^{q} v_+^{q+1}\, h \Big)\,dx, \\[4pt]
    \langle R_2(u,v), h\rangle
    &:= \int_{\R^n}\Big( \nabla v\cdot\nabla h - v_+^{p}\, h
        - \beta\, u_+^{q+1} v_+^{q}\, h \Big)\,dx,
\end{aligned}
\end{equation}
for $h\in \dot H^1(\R^n)$. Both are bounded linear functionals on
$\dot H^1(\R^n)$ by H\"older's inequality and the Sobolev embedding, and
$(u,v)$ solves \eqref{eq:positive-system} if and only if
$$
    R_1(u,v) = 0, \qquad R_2(u,v) = 0
$$
as elements of the dual space. Equivalently, writing
$w_i := (-\Delta)^{-1}R_i(u,v)$ for the Riesz representative, one has
$\langle R_i(u,v),h\rangle = \langle w_i, h\rangle_{H}$ for every $h$, so that the
two formulations agree.

The corrections will be sought in $X_k$, but the solutions obtained are weak solutions in the
whole of $\dot H^1(\R^n)$, by the following observation.

\begin{remark}\label{rem:symmetric-criticality}
Let $u,v\in X_k$. Then $w_i=(-\Delta)^{-1}R_i(u,v)\in X_k$ for $i=1,2$. Indeed,
$$
w_1=u-(-\Delta)^{-1}f_1,\qquad w_2=v-(-\Delta)^{-1}f_2,\qquad
f_1:=u_+^p+\beta u_+^qv_+^{q+1},\quad f_2:=v_+^p+\beta u_+^{q+1}v_+^{q}.
$$
The reflections and the rotation in the definition of $X_k$ commute with $-\Delta$ and leave
$f_1,f_2$ invariant. For the Kelvin transform, the following identity holds distributionally:
$$
-\Delta(Kw)(x)=|x|^{-(n+2)}\,(-\Delta w)\Big(\frac{x}{|x|^2}\Big).
$$
By uniqueness of the Riesz representative, the function $(-\Delta)^{-1}f$ is
Kelvin invariant whenever
$f(x)=|x|^{-(n+2)}f(x/|x|^2)$. Since $u=Ku$, $v=Kv$ and $(n-2)p=(n-2)(2q+1)=n+2$,
$$
u_+^p(x)=|x|^{-(n+2)}\,u_+^p\Big(\frac{x}{|x|^2}\Big),
\qquad
\big(u_+^qv_+^{q+1}\big)(x)=|x|^{-(n+2)}\,\big(u_+^qv_+^{q+1}\big)\Big(\frac{x}{|x|^2}\Big),
$$
and similarly for the terms of $f_2$. Hence $w_1,w_2\in X_k$. As a consequence, if $R_i(u,v)$
vanishes on $X_k$, then, testing with $w_i\in X_k$,
$$
\|w_i\|^2=\big\langle R_i(u,v),w_i\big\rangle=0,
$$
that is, $R_i(u,v)=0$ in the dual of $\dot H^1(\R^n)$. In other words, a pair
$(u,v)\in X_k\times X_k$ which solves \eqref{eq:positive-system} when tested against $X_k$ is a
weak solution in the whole of $\dot H^1(\R^n)$.
\end{remark}

For fixed parameters $(a,b,\delta)$, the \emph{infinite-dimensional part} of the
reduction consists in finding $(\phi,\psi)\in\mathcal H_{1}\times\mathcal
H_{2}$ such that
\begin{equation}\label{eq:correction-eqs}
\begin{cases}
    \langle R_1(aU+\phi,\,bV+\psi),\, h\rangle = 0
        &\text{for all } h \in \mathcal H_{1}, \\[4pt]
    \langle R_2(aU+\phi,\,bV+\psi),\, g\rangle = 0
        &\text{for all } g \in \mathcal H_{2}.
\end{cases}
\end{equation}
Identifying the dual space with $\dot H^1(\R^n)$ through the Riesz isomorphism,
we may equivalently write, in terms of the projections,
$$
    P_1^\perp R_1(aU+\phi, bV+\psi) = 0,
    \qquad
    P_2^\perp R_2(aU+\phi, bV+\psi) = 0 .
$$
This problem is variational and is solved in Section~\ref{sec:infinite-dim}. 
It remains to impose the components of the residuals along the finite-dimensional
spaces $\mathcal F_1 = \operatorname{span}\{U\}$ and $\mathcal F_2 =
\operatorname{span}\{V,Z_\delta\}$, that is, the \emph{bifurcation equations}
\begin{equation}\label{eq:bifurcation-eqs}
\begin{cases}
    \big\langle R_1(aU+\phi,\,bV+\psi),\, U \big\rangle = 0, \\[4pt]
    \big\langle R_2(aU+\phi,\,bV+\psi),\, V \big\rangle = 0, \\[4pt]
    \big\langle R_2(aU+\phi,\,bV+\psi),\, Z_\delta \big\rangle = 0,
\end{cases}
\end{equation}
equivalently $P_1 R_1 = 0$ and $P_2 R_2 = 0$. This is a system of three equations
in the three unknowns $(a,b,\delta)$, solved in Section~\ref{sec:finite-dim}.

\section{The infinite-dimensional reduction via minimization}\label{sec:infinite-dim}

\subsection{Preliminary error estimates}
Before performing the infinite-dimensional reduction, we estimate how accurately the ansatz $(aU,bV)$ solves the system. Let us denote
\begin{equation}\label{eq:errors}
\begin{cases}
    E_1 :=\Delta (aU) + (aU)^{p} + \beta\, (aU)^{q} (bV)^{q+1}, \\[4pt]
    E_2 := \Delta (bV) + (bV)^{p} + \beta\, (aU)^{q+1} (bV)^{q}.
\end{cases}
\end{equation}

  \begin{proposition}\label{prop:error}

Let $\varepsilon$ and $(a,b,\delta)$ be given by
\eqref{eq:eps-beta} and \eqref{eq:rescaling}. Assume that $d$ ranges in a fixed compact subset of $(0,\infty)$
and that $t,s$ remain bounded. Set
$$
\eta_n(\varepsilon):=
\begin{cases}
\varepsilon^3, & n=5,\\
\varepsilon^4|\ln\varepsilon|^{2/3}, & n=6,\\
\varepsilon^{(n+2)/2}, & n\ge7.
\end{cases}
$$
Then, as $\beta\to0^-$,
$$
\big\|U^qV^{q+1}\big\|_{(2^*)'}
+\big\|U^{q+1}V^q\big\|_{(2^*)'}
\lesssim\varepsilon,
\qquad
\big\|\Delta V+V^p\big\|_{(2^*)'}
\lesssim\eta_n(\varepsilon).
$$
Consequently,
$$
\|E_1\|_{(2^*)'}
\lesssim
\varepsilon^{n-2}
+\frac{\varepsilon^{(n-2)/2}}{|\ln\varepsilon|},
$$
and
$$
\|E_2\|_{(2^*)'}
\lesssim
\varepsilon^{n-2}
+\eta_n(\varepsilon)
+\frac{\varepsilon^{(n-2)/2}}{|\ln\varepsilon|}.
$$
In particular,
$$
\|E_1\|_{(2^*)'}+\|E_2\|_{(2^*)'}
=
O\!\left(
\frac{\varepsilon^{(n-2)/2}}{|\ln\varepsilon|}
\right).
$$
All estimates are uniform on the specified parameter set.
\end{proposition}

\begin{proof}
Set
$$
r_*:=(2^*)'=\frac{2n}{n+2}.
$$
Throughout the proof, the implicit constants may depend on $n$, $k$,
the compact set in which $d$ ranges, and the bounds on $t,s$,
but are independent of $\varepsilon$.

Let $\mathcal R$ denote the rotation through the angle $2\pi/k$
in the $(x_1,x_2)$-plane, and set
\begin{equation}\label{eq:Omega}
\Omega_1:=
\big\{x\in\R^n:\ |x-\xi_1|<|x-\xi_j|
\text{ for every }j\ne1\big\},
\qquad
\Omega_j:=\mathcal R^{j-1}\Omega_1,
\qquad
W:=\sum_{j=2}^kU_j.
\end{equation}
These regions partition $\R^n$ up to a set of measure zero.
On $\Omega_1$ we have $U_j\le U_1$ for every $j$, and hence
$V\le kU_1$.

\medskip
\noindent
\emph{Step 1: the coupling terms.}
Set
$$
\kappa_1:=(q+1)r_*=\frac{2n^2}{(n-2)(n+2)},
\qquad
\kappa_2:=qr_*=\frac{4n}{(n-2)(n+2)}.
$$
Since $\kappa_1+\kappa_2=2^*$, we have
$$
\|U^qV^{q+1}\|_{r_*}^{r_*}
=\int_{\R^n}U^{\kappa_2}V^{\kappa_1},
\qquad
\|U^{q+1}V^q\|_{r_*}^{r_*}
=\int_{\R^n}U^{\kappa_1}V^{\kappa_2}.
$$
The integrands are invariant under $\mathcal R$. Therefore,
for
$(\alpha_1,\alpha_2)\in
\{(\kappa_2,\kappa_1),(\kappa_1,\kappa_2)\}$,
$$
\int_{\R^n}U^{\alpha_1}V^{\alpha_2}
=k\int_{\Omega_1}U^{\alpha_1}V^{\alpha_2}
\le k^{1+\alpha_2}
\int_{\R^n}U^{\alpha_1}U_1^{\alpha_2}.
$$
For the first integral, using $U\le C_n$ and the change of
variables $x=\xi_1+\delta y$, we obtain
$$
\int_{\R^n}U^{\kappa_2}U_1^{\kappa_1}
\lesssim
\delta^{\,n-\kappa_1(n-2)/2}
\int_{\R^n}U^{\kappa_1}
\lesssim\delta^{r_*}.
$$
Indeed,
$$
n-\frac{\kappa_1(n-2)}2=r_*,
\qquad
\kappa_1(n-2)=\frac{2n^2}{n+2}>n,
$$
so that $\int_{\R^n}U^{\kappa_1}<\infty$.

For the second integral, the pointwise estimate
$$
U_1(x)\le C_n\delta^{(n-2)/2}|x-\xi_1|^{2-n}
$$
gives
$$
\int_{\R^n}U^{\kappa_1}U_1^{\kappa_2}
\lesssim
\delta^{r_*}
\int_{\R^n}
\frac{dx}
{(1+|x|^2)^{n^2/(n+2)}
 |x-\xi_1|^{4n/(n+2)}}
\lesssim\delta^{r_*}.
$$
The last integral is uniformly bounded for $|\xi_1|\le1$:
on $|x-\xi_1|<1$ this follows from
$4n/(n+2)<n$, while on the complement it follows from
$2n^2/(n+2)>n$.

Taking the $r_*$-th root, we conclude that
\begin{equation}\label{eq:coupling-delta-bound}
\|U^qV^{q+1}\|_{r_*}
+\|U^{q+1}V^q\|_{r_*}
\lesssim\delta.
\end{equation}

\medskip
\noindent
\emph{Step 2: the interaction error.}
Since $-\Delta U_j=U_j^p$, we have
$$
N:=\Delta V+V^p=V^p-\sum_{j=1}^kU_j^p.
$$
On $\Omega_1$, using $W\le(k-1)U_1$ and $p>1$, we obtain
$$
0\le N
\le (U_1+W)^p-U_1^p
\le p(U_1+W)^{p-1}W
\lesssim U_1^{p-1}W.
$$
By rotational symmetry,
$$
\|N\|_{r_*}
\lesssim
\|U_1^{p-1}W\|_{L^{r_*}(\Omega_1)}.
$$
Write $r:=|x-\xi_1|$. For $x\in\Omega_1$ and $j\ne1$,
$$
|x-\xi_j|\ge r,
\qquad
|x-\xi_j|\ge\frac12|\xi_1-\xi_j|.
$$
Since $k$ is fixed and $\rho\to1$, the distances
$|\xi_1-\xi_j|$, $j\ne1$, are bounded away from zero
for $\delta$ sufficiently small. Consequently,
$$
W(x)\lesssim
\delta^{(n-2)/2}\min\{1,r^{2-n}\}
\qquad\text{on }\Omega_1.
$$
In the region $r<1$, we use $W\lesssim\delta^{(n-2)/2}$.
Setting
$$
m:=(p-1)r_*=\frac{8n}{(n-2)(n+2)}
$$
and changing variables $x=\xi_1+\delta y$, we get
$$
\int_{\Omega_1\cap\{r<1\}}(U_1^{p-1}W)^{r_*}
\lesssim
\delta^{(n-2)r_*/2}
\delta^{\,n-m(n-2)/2}
\int_{|y|<\delta^{-1}}U(y)^m\,dy.
$$
Here
$$
n-\frac{m(n-2)}2=\frac{n(n-2)}{n+2},
\qquad
m(n-2)-n=\frac{n(6-n)}{n+2}.
$$
It follows that
$$
\int_{|y|<\delta^{-1}}U(y)^m\,dy
\lesssim
\begin{cases}
1, & n=5,\\
|\ln\delta|, & n=6,\\
\delta^{-n(n-6)/(n+2)}, & n\ge7.
\end{cases}
$$
Taking the $r_*$-th root yields
$$
\|U_1^{p-1}W\|_{L^{r_*}(\Omega_1\cap\{r<1\})}
\lesssim
\begin{cases}
\delta^3, & n=5,\\
\delta^4|\ln\delta|^{2/3}, & n=6,\\
\delta^{(n+2)/2}, & n\ge7.
\end{cases}
$$
In the region $r\ge1$, we use
$$
U_1^{p-1}\lesssim\delta^2r^{-4},
\qquad
W\lesssim\delta^{(n-2)/2}r^{2-n}.
$$
Thus
$$
\int_{\Omega_1\cap\{r\ge1\}}(U_1^{p-1}W)^{r_*}
\lesssim
\delta^n\int_{\{|x-\xi_1|\ge1\}}
|x-\xi_1|^{-2n}\,dx
\lesssim\delta^n.
$$
This region therefore contributes at most
$O(\delta^{(n+2)/2})$ to the norm. Combining the two regions,
we obtain
\begin{equation}\label{eq:defect-delta-bound}
\|\Delta V+V^p\|_{r_*}
\lesssim
\begin{cases}
\delta^3, & n=5,\\
\delta^4|\ln\delta|^{2/3}, & n=6,\\
\delta^{(n+2)/2}, & n\ge7.
\end{cases}
\end{equation}

\medskip
\noindent
\emph{Step 3: the errors of the full ansatz.}
For $\varepsilon$ sufficiently small, $a,b\in[1/2,2]$.
The definitions of $E_1,E_2$ give the exact identities
$$
E_1=(a^p-a)U^p
+\beta a^qb^{q+1}U^qV^{q+1},
$$
$$
E_2=b(\Delta V+V^p)+(b^p-b)V^p
+\beta a^{q+1}b^qU^{q+1}V^q.
$$
Moreover,
$$
|a^p-a|\lesssim|a-1|,
\qquad
|b^p-b|\lesssim|b-1|,
$$
and
$$
\|U^p\|_{r_*}=\|U\|_{2^*}^p,
\qquad
\|V^p\|_{r_*}=\|V\|_{2^*}^p
\le k^p\|U\|_{2^*}^p.
$$
Hence, by \eqref{eq:coupling-delta-bound},
$$
\|E_1\|_{r_*}\lesssim |a-1|+|\beta|\delta
\qquad\text{and}\qquad
\|E_2\|_{r_*}
\lesssim |b-1|+\|\Delta V+V^p\|_{r_*}
+|\beta|\delta.
$$
Finally, using \eqref{eq:rescaling} and \eqref{eq:eps-beta},
$$
|a-1|+|b-1|\lesssim\varepsilon^{n-2},
\qquad
|\beta|\delta
=d\,\frac{\varepsilon^{(n-2)/2}}{|\ln\varepsilon|}.
$$
Since $d$ ranges in a fixed compact subset of $(0,\infty)$,
we have $\delta\asymp\varepsilon$ and
$|\ln\delta|\asymp|\ln\varepsilon|$. Therefore
\eqref{eq:coupling-delta-bound} and
\eqref{eq:defect-delta-bound} imply
$$
\|U^qV^{q+1}\|_{r_*}
+\|U^{q+1}V^q\|_{r_*}\lesssim\varepsilon,
\qquad
\|\Delta V+V^p\|_{r_*}\lesssim\eta_n(\varepsilon).
$$
Consequently,
$$
\|E_1\|_{r_*}
\lesssim
\varepsilon^{n-2}
+\frac{\varepsilon^{(n-2)/2}}{|\ln\varepsilon|},
$$
$$
\|E_2\|_{r_*}
\lesssim
\varepsilon^{n-2}+\eta_n(\varepsilon)
+\frac{\varepsilon^{(n-2)/2}}{|\ln\varepsilon|}.
$$
Since, for every $n\ge5$,
$$
\varepsilon^{n-2}+\eta_n(\varepsilon)
=o\!\left(
\frac{\varepsilon^{(n-2)/2}}{|\ln\varepsilon|}
\right),
$$
the final estimate follows.

\end{proof}

\subsection{The quadratic form}

Set
$$
Q_1(\phi):=\|\phi\|^2-p\int_{\R^n}U^{p-1}\phi^2,
\qquad
Q_2(\psi):=\|\psi\|^2-p\int_{\R^n}V^{p-1}\psi^2,
$$
and
$$
Q(\phi,\psi):=\frac12Q_1(\phi)+\frac12Q_2(\psi).
$$
We claim that there exists $\delta_0>0$ such that, for every $\delta\in(0,\delta_0)$,
\begin{equation}\label{eq:coercivity}
Q(\phi,\psi)\ge
\gamma\big(\|\phi\|^2+\|\psi\|^2\big)
\qquad
\text{for all }(\phi,\psi)\in\mathcal H_1\times\mathcal H_2,
\qquad
\gamma>0,
\end{equation}
$\gamma$ independent of $\delta$.

We first recall the spectral properties of the linearized
operator around the standard bubble in the full space
$\dot H^1(\R^n)$, that can be found in \cite{Rey1990,BianchiEgnell1991}. Consider the weighted eigenvalue problem
$$
-\Delta\phi-pU^{p-1}\phi=\lambda U^{p-1}\phi,
\qquad \phi\in\dot H^1(\R^n).
$$
Its first eigenvalue is $1-p<0$, with eigenspace
$\operatorname{span}\{U\}$, while its kernel is
$
\operatorname{span}\{Z^0,Z^1,\ldots,Z^n\},
$
where the functions $Z^i$ are defined in \eqref{Zi}. Moreover, the spectral gap above this kernel gives
$$
\|\phi\|^2-p\int_{\R^n}U^{p-1}\phi^2
\ge \frac4{n+4}\|\phi\|^2,
$$
whenever $\phi$ is orthogonal, with respect to
$\langle\cdot,\cdot\rangle_H$, to $U$ and to every $Z^i$;
see \cite[Appendix~D]{Rey1990}.

We now restrict the quadratic form to $X_k$.
The symmetries defining $X_k$ imply
$$
\langle\phi,Z^i\rangle_H=0
\qquad
\text{for every }\phi\in X_k,\quad i=0,\ldots,n.
$$
Indeed, the Kelvin transform is an isometry of $\dot H^1(\R^n)$
and satisfies $KZ^0=-Z^0$, whereas $K\phi=\phi$.
The orthogonality to $Z^i$, $i=2,\ldots,n$, follows from
the evenness of $\phi$ in $x_i$ and the oddness of $Z^i$
in that variable. Finally, averaging over the rotations
through multiples of $2\pi/k$ gives the orthogonality to
$Z^1$, since $k\ge2$ and the rotational average of $Z^1$
vanishes.
Consequently,
\begin{equation}\label{q1prev}
\lambda_*:=
\inf_{\substack{
\phi\in X_k,\ \phi\ne0\\
\langle\phi,U\rangle_H=0}}
\frac{\|\phi\|^2-p\int_{\R^n}U^{p-1}\phi^2}
{\|\phi\|^2}
\ge \frac4{n+4}.
\end{equation}
In particular, with $\gamma_1:=4/(n+4)$,
\begin{equation}\label{q1}
Q_1(\phi)\ge\gamma_1\|\phi\|^2
\qquad\text{for all }\phi\in\mathcal H_1.
\end{equation}

\begin{lemma} \label{lem:spg}
There exist
$\delta_0=\delta_0(n,k)>0$ and $\mu_0=\mu_0(n,k)>0$ such that,
for every $\delta\in(0,\delta_0)$,
\begin{equation}\label{mu*}
\mu_*(\delta):=
\inf_{\substack{
\psi\in X_k,\ \psi\ne0\\
\langle\psi,V\rangle_H=0,\ 
\langle\psi,Z_\delta\rangle_H=0}}
\frac{\|\psi\|^2-p\int_{\R^n}V^{p-1}\psi^2}
{\|\psi\|^2}
\ge\mu_0.
\end{equation}
\end{lemma}

\begin{proof}
We prove the lemma by contradiction and blow-up: by the
symmetries of $X_k$, a function in $\mathcal H_2$ looks, near each
$\xi_j$ and after rescaling, like a function orthogonal to $U$ and to
all the $Z^i$, for which the single-bubble bound above applies.

Suppose there are $\delta_m\to0^+$ and $\psi_m\in\mathcal H_2$,
with $\delta=\delta_m$, such that $\|\psi_m\|=1$ and $Q_2(\psi_m)\le 1/m$,
so that 
\begin{equation}\label{lowBound}
p\int_{\R^n}V^{p-1}\psi_m^2\ge1-1/m.
\end{equation}
Let $\mathcal R$ be the rotation through the angle
$2\pi/k$ in the $(x_1,x_2)$-plane. Since 
$$U_j=U_1\circ\mathcal R^{1-j},\quad Z_j=Z_1\circ\mathcal R^{1-j},\quad \psi_m\circ\mathcal R=\psi_m,
$$
then
$$\langle\psi_m,U_1\rangle_H=\tfrac1k\langle\psi_m,V\rangle_H=0,
\qquad
\langle\psi_m,Z_1\rangle_H=\tfrac1k\langle\psi_m,Z_\delta\rangle_H=0 .
$$
Fix $R>0$. For $m$ large the balls $\mathcal B_j:=B(\xi_j,R\delta)$ are
pairwise disjoint and permuted by $\mathcal R$. Hence, 
$$
\int_{\mathcal B_1}|\nabla\psi_m|^2\le\frac1k,
$$
and 
\begin{equation}\label{splitV}
p\int_{\R^n}V^{p-1}\psi_m^2
=k\,p\int_{\mathcal B_1}V^{p-1}\psi_m^2+p\int_{\R^n\setminus\bigcup_j\mathcal B_j}V^{p-1}\psi_m^2.
\end{equation}
Since $V^{2^*}\le k^{2^*-1}\sum_jU_j^{2^*}$ and $|x-\xi_j|\ge R\delta$
on $\R^n\setminus\bigcup_j\mathcal B_j$, H\"older's and Sobolev's inequalities imply
$$\int_{\R^n\setminus\bigcup_j\mathcal B_j}V^{p-1}\psi_m^2\leq C\big(\displaystyle\int_{\{|y|\ge R\}}U^{2^*}\big)^{2/n},$$ 
uniformly in $m$.

Let $T_mf(y):=\delta^{(n-2)/2}f(\xi_1+\delta y)$, an
isometry of $\dot H^1(\R^n)$. By scaling, $T_mU_1=U$,
$T_mZ_1=Z^0-\frac\delta\rho Z^1$, and
$$T_mV=U+\sum_{j\ge2}U\big(\cdot+\frac{\xi_1-\xi_j}{\delta}\big)\to U\quad \mbox{
locally uniformly}.$$ 
Up to a subsequence,
$\widetilde\psi_m:=T_m\psi_m\rightharpoonup\widetilde\psi$ weakly in
$\dot H^1(\R^n)$ and strongly in $L^2_{\rm loc}(\R^n)$. By weak lower
semicontinuity,
$$\int_{B(0,R)}|\nabla\widetilde\psi|^2
\le\liminf_m\int_{\mathcal B_1}|\nabla\psi_m|^2\le\frac1k$$ 
for every
$R$, so $\|\widetilde\psi\|^2\le\frac1k.$ Moreover,
$$\int_{\mathcal B_1}V^{p-1}\psi_m^2
=\int_{B(0,R)}(T_mV)^{p-1}\widetilde\psi_m^2
\to\int_{B(0,R)}U^{p-1}\widetilde\psi^2,$$ so, using \eqref{lowBound}, and letting first $m\to\infty$
and then $R\to\infty$ in \eqref{splitV}, we conclude
\begin{equation}\label{contrad}
p\int_{\R^n}U^{p-1}\widetilde\psi^2\ge\frac1k.
\end{equation}
Since $T_m$ is an isometry,
$$
\langle\widetilde\psi_m,U\rangle_H=\langle\psi_m,U_1\rangle_H
=\tfrac1k\langle\psi_m,V\rangle_H=0,
\qquad
\langle\widetilde\psi_m,Z^0\rangle_H
=\tfrac\delta\rho\langle\widetilde\psi_m,Z^1\rangle_H=O(\delta);
$$
hence $\widetilde\psi$ is orthogonal to $U$ and $Z^0$. The evenness of
$\psi_m$ in $x_2,\dots,x_n$ passes to $\widetilde\psi$, which is
therefore orthogonal to $Z^2,\dots,Z^n$. Finally, the Kelvin invariance
reads 
$$\widetilde\psi_m(y)=|\xi_1+\delta y|^{2-n}
\widetilde\psi_m(\Phi_m(y)),\quad \mbox{ with }\quad \Phi_m(y):=\delta^{-1}\big(\iota(\xi_1+\delta y)-\xi_1\big),\quad \iota(x):=x/|x|^2.$$ 
Since $\rho^2+\delta^2=1$ gives
$\iota(\xi_1)-\xi_1=\delta^2\xi_1/\rho^2$, the maps $\Phi_m$ converge
in $C^1_{\rm loc}$ to the reflection $y_1\mapsto-y_1$. Testing against
$\varphi\in C_c^\infty(\R^n)$ and changing variables $y=\Phi_m(z)$, we
conclude that $\widetilde\psi$ is even in $y_1$, hence orthogonal to
$Z^1$.

By \eqref{q1prev},
$$0=\frac1k-\frac1k\geq \|\widetilde\psi\|^2-p\int_{\R^n}U^{p-1}\widetilde\psi^2
\ge\frac4{n+4}\|\widetilde\psi\|^2.$$ 
Thus $\widetilde\psi=0$, which contradicts \eqref{contrad}.
\end{proof}

Lemma \ref{lem:spg} immediately yields, with $\gamma_2:=\mu_0$,
\begin{equation}\label{q2}
Q_2(\psi)\ge\gamma_2\|\psi\|^2
\qquad\text{for all }\psi\in\mathcal H_2,
\end{equation}
uniformly for $\delta\in(0,\delta_0)$.
Together with \eqref{q1}, this proves \eqref{eq:coercivity} with $\gamma:=\frac12\min\{\gamma_1,\gamma_2\}$.

\subsection{Existence of a minimum}
Let us consider the energy
$$
\begin{aligned}
J_\beta(u,v)
&=\frac12\int_{\R^n}(|\nabla u|^2+|\nabla v|^2)\,dx
-\frac1{p+1}\int_{\R^n}(u_+^{p+1}+v_+^{p+1})\,dx-\frac{2\beta}{p+1}\int_{\R^n}u_+^{(p+1)/2}v_+^{(p+1)/2}\,dx.
\end{aligned}
$$
For short, let us call
\begin{equation}\label{defF}
F(s):=\frac 1{p+1}(s _+)^{p+1},\quad f(s)=F'(s),\quad f'(s)=F''(s),
\end{equation}
and
\begin{equation}\label{defG}
G(s,t):=\frac 2{p+1}(s _+)^{\frac{p+1}2}(t _+)^{\frac{p+1}2}.
\end{equation}
Returning to the energy, set
$$
I_{a,b,\delta}(\phi,\psi)
:=J_\beta(aU+\phi,bV+\psi)-J_\beta(aU,bV),
\qquad
(\phi,\psi)\in\mathcal H_1\times\mathcal H_2.
$$
For brevity, we write $I=I_{a,b,\delta}$, and throughout this subsection, the parameters are chosen as in
\eqref{eq:rescaling}, with $\varepsilon$ defined by
\eqref{eq:eps-beta}. 
Since $\langle R_1(u,v),\cdot\rangle$ and $\langle R_2(u,v),\cdot\rangle$ (defined in \eqref{eq:residuals}) are
precisely the partial derivatives $\partial_u J_\beta(u,v)$ and $\partial_v
J_\beta(u,v)$, the equations \eqref{eq:correction-eqs} are exactly the
critical point equations for $I_{a,b,\delta}$ restricted to
$\mathcal H_1\times\mathcal H_2$.

We will prove the existence of a critical point of $I$ by following the local minimization argument of
\cite[Proposition~2]{BahriLiRey1995} and
\cite[Proposition~4]{Rey1990}, which goes back to
\cite[Lemma~4.1]{Bahri1989}.
For $r>0$, write
$$
B_\delta(0,r):=
\big\{(\phi,\psi)\in\mathcal H_1\times\mathcal H_2:
\|\phi\|+\|\psi\|\le r\big\}.
$$
The subindex $\delta$ points out the dependence of $\mathcal H_2$ on this parameter.
\begin{lemma}
\label{lem:convexity}
Fix a compact set $D\subset(0,\infty)$ and $T>0$.
There exist $\beta_1<0$ and $r_0\in(0,1)$ such that,
for every $\beta\in[\beta_1,0)$ and every choice of
parameters as in \eqref{eq:rescaling} with
$d\in D$ and $|t|+|s|\le T$, the following holds:
for all
$(\phi_i,\psi_i)\in B_\delta(0,r_0)$, $i=1,2$,
\begin{equation}\label{eq:strict-convexity}
\begin{aligned}
&\Big\langle
\nabla I(\phi_1,\psi_1)-\nabla I(\phi_2,\psi_2),
(\phi_1-\phi_2,\psi_1-\psi_2)
\Big\rangle_H\ge
\gamma\big(
\|\phi_1-\phi_2\|^2+\|\psi_1-\psi_2\|^2
\big),
\end{aligned}
\end{equation}
where $\gamma$ is as in \eqref{eq:coercivity}.
Here $\nabla I$ denotes the gradient of the restriction
of $I$ to $\mathcal H_1\times\mathcal H_2$, equipped with
the product Hilbert inner product.
\end{lemma}

\begin{proof}
Set
$$
h:=\phi_1-\phi_2,\qquad g:=\psi_1-\psi_2,
\qquad
u_i:=aU+\phi_i,\qquad v_i:=bV+\psi_i.
$$
Since $I$ is of class $C^1$, its first derivative gives
$$
\begin{aligned}
\big\langle
\nabla I(\phi_1,\psi_1)-\nabla I(\phi_2,\psi_2),
(h,g)\big\rangle_H
&=\|h\|^2-pa^{p-1}\int_{\R^n}U^{p-1}h^2
+\|g\|^2-pb^{p-1}\int_{\R^n}V^{p-1}g^2\\
&\quad-\mathcal T_\phi-\mathcal T_\psi
-\beta\int_{\R^n}
\big\langle\nabla G(u_1,v_1)-\nabla G(u_2,v_2),
(h,g)\big\rangle_{\R^2}\,dx,
\end{aligned}
$$
where
$$
\mathcal T_\phi:=
\int_{\R^n}
\big[f(aU+\phi_1)-f(aU+\phi_2)-f'(aU)h\big]h\,dx,
$$
$$
\mathcal T_\psi:=
\int_{\R^n}
\big[f(bV+\psi_1)-f(bV+\psi_2)-f'(bV)g\big]g\,dx,
$$
and $f$ and $G$ are defined in \eqref{defF} and \eqref{defG}.
\smallskip 

\noindent
\emph{The quadratic terms.}
Since $(h,g)\in\mathcal H_1\times\mathcal H_2$,
\eqref{eq:coercivity} gives
$$
2Q(h,g)\ge2\gamma\big(\|h\|^2+\|g\|^2\big).
$$
Furthermore,
$$
|a^{p-1}-1|+|b^{p-1}-1|
\lesssim |a-1|+|b-1|
\lesssim\varepsilon^{n-2}.
$$
Using H\"older's and Sobolev's inequalities, we obtain
$$
\begin{aligned}
&\|h\|^2-pa^{p-1}\int_{\R^n}U^{p-1}h^2
+\|g\|^2-pb^{p-1}\int_{\R^n}V^{p-1}g^2\ge
\big(2\gamma-C\varepsilon^{n-2}\big)
\big(\|h\|^2+\|g\|^2\big).
\end{aligned}
$$
\smallskip

\noindent
\emph{The scalar remainders.}
By the fundamental theorem of calculus,
$$
\mathcal T_\phi
=\int_{\R^n}\int_0^1
\big[f'(aU+w_\tau)-f'(aU)\big]h^2\,d\tau\,dx,
\qquad
w_\tau:=\phi_2+\tau h,
$$
and, since $\|\phi_i\|\le r_0$, we have
$\|w_\tau\|\le r_0$.

If $n\ge6$, then
$$
|f'(aU+w_\tau)-f'(aU)|
\lesssim|w_\tau|^{p-1}.
$$
H\"older's inequality gives
$$
\int_{\R^n}|w_\tau|^{p-1}h^2
\le\|w_\tau\|_{2^*}^{p-1}\|h\|_{2^*}^2
\lesssim r_0^{p-1}\|h\|^2.
$$

If $n=5$, then
$$
|f'(aU+w_\tau)-f'(aU)|
\lesssim U^{p-2}|w_\tau|+|w_\tau|^{p-1}.
$$
In addition to the preceding estimate, we use
$$
\int_{\R^n}U^{p-2}|w_\tau|h^2
\le
\|U\|_{2^*}^{p-2}
\|w_\tau\|_{2^*}\|h\|_{2^*}^2
\lesssim r_0\|h\|^2.
$$
Thus, with $m:=\min\{p-1,1\}$ and $r_0<1$,
$$
|\mathcal T_\phi|\le Cr_0^m\|h\|^2.
$$
The same argument, using $\|V\|_{2^*}\le k\|U\|_{2^*}$,
gives
$$
|\mathcal T_\psi|\le Cr_0^m\|g\|^2.
$$

\medskip
\noindent
\emph{The coupling term.}
Recall that
$$
\partial_sG(s,t)=s_+^qt_+^{q+1},
\qquad
\partial_tG(s,t)=s_+^{q+1}t_+^q.
$$
Pointwise, we have the exact decomposition
$$
\begin{aligned}
\big\langle\nabla G(u_1,v_1)-\nabla G(u_2,v_2),
(h,g)\big\rangle_{\R^2}
&=
\big((u_1)_+^q-(u_2)_+^q\big)(v_1)_+^{q+1}h+(u_2)_+^q
\big((v_1)_+^{q+1}-(v_2)_+^{q+1}\big)h\\
&\quad+\big((v_1)_+^q-(v_2)_+^q\big)(u_2)_+^{q+1}g+(v_1)_+^q
\big((u_1)_+^{q+1}-(u_2)_+^{q+1}\big)g.
\end{aligned}
$$
The first and third terms are nonnegative, since
$h=u_1-u_2$, $g=v_1-v_2$, and $s\mapsto s_+^q$
is nondecreasing.
For the other two terms, we use
$$
|x_+^{q+1}-y_+^{q+1}|
\lesssim (|x|^q+|y|^q)|x-y|.
$$
Hence 
$$
|(u_2)_+^q
\big((v_1)_+^{q+1}-(v_2)_+^{q+1}\big)h+(v_1)_+^q
\big((u_1)_+^{q+1}-(u_2)_+^{q+1}\big)g|\leq C\big(|u_1|^q+|u_2|^q\big)
\big(|v_1|^q+|v_2|^q\big)|h||g|.
$$
Since $2q+2=2^*$, H\"older's inequality gives
$$
\begin{aligned}
\int_{\R^n}
\big(|u_1|^q+|u_2|^q\big)
\big(|v_1|^q+|v_2|^q\big)|h||g|
&\le
\big(\|u_1\|_{2^*}^q+\|u_2\|_{2^*}^q\big)
\big(\|v_1\|_{2^*}^q+\|v_2\|_{2^*}^q\big)
\|h\|_{2^*}\|g\|_{2^*}\\
&\lesssim
\|h\|\|g\|
\lesssim\|h\|^2+\|g\|^2.
\end{aligned}
$$
Here the norms of $u_i,v_i$ are uniformly bounded because
$a,b\le2$, $\|V\|_{2^*}\le k\|U\|_{2^*}$ and
$\|\phi_i\|+\|\psi_i\|\le r_0<1$.
Since $-\beta=|\beta|>0$, we conclude that
$$
\begin{aligned}
&-\beta\int_{\R^n}
\big\langle\nabla G(u_1,v_1)-\nabla G(u_2,v_2),
(h,g)\big\rangle_{\R^2}\,dx
\ge
-C|\beta|\big(\|h\|^2+\|g\|^2\big).
\end{aligned}
$$
Combining the preceding estimates yields
$$
\begin{aligned}
&\big\langle
\nabla I(\phi_1,\psi_1)-\nabla I(\phi_2,\psi_2),
(h,g)\big\rangle_H\ge
\big(2\gamma-C\varepsilon^{n-2}
-Cr_0^m-C|\beta|\big)
\big(\|h\|^2+\|g\|^2\big).
\end{aligned}
$$
Choose $r_0\in(0,1)$ so that $Cr_0^m\le\gamma/2$.
Since $\varepsilon\to0$ as $\beta\to0^-$, we may then
choose $\beta_1<0$ sufficiently close to zero so that
$$
C\varepsilon^{n-2}+C|\beta|\le\gamma/2
$$
and the coercivity estimate applies throughout the
parameter set. This proves \eqref{eq:strict-convexity}.
\end{proof}

\begin{proposition}
\label{lem:minimizer}
Fix a compact set $D\subset(0,\infty)$ and $T>0$.
Let $\varepsilon$ be defined by \eqref{eq:eps-beta},
and choose the parameters as in \eqref{eq:rescaling},
with
$$
d\in D,\qquad |t|+|s|\le T.
$$
There exist constants $R>0$ and $\beta_0<0$ such that,
for every $\beta\in[\beta_0,0)$, the functional
$I_{a,b,\delta}$ has a unique minimum point
$$
(\phi_{a,b,\delta},\psi_{a,b,\delta})
\in B_\delta\left(0,R\frac{\varepsilon^{(n-2)/2}}{|\ln\varepsilon|}\right).
$$
This point lies in the interior of the ball and solves
\eqref{eq:correction-eqs}. Moreover, it is the unique
critical point of $I_{a,b,\delta}$ in $B_\delta(0,r_0)$,
where $r_0$ is given by Lemma~\ref{lem:convexity}, and
$$
\|\phi_{a,b,\delta}\|+\|\psi_{a,b,\delta}\|
\lesssim\frac{\varepsilon^{(n-2)/2}}{|\ln\varepsilon|}.
$$
All conclusions are uniform on the specified parameter set.
\end{proposition}

\begin{proof}
Write $I=I_{a,b,\delta}$ and equip
$\mathcal H_1\times\mathcal H_2$ with the product Hilbert norm
$$
\|z\|_H:=
\big(\|\phi\|^2+\|\psi\|^2\big)^{1/2},
\qquad z=(\phi,\psi).
$$
Thus
$$
\|z\|_H\le\|\phi\|+\|\psi\|
\le\sqrt2\,\|z\|_H.
$$
By Proposition~\ref{prop:error}, H\"older's inequality
and the Sobolev inequality, there exists a uniform
constant $C_0>0$ such that
$$
\begin{aligned}
|\langle\nabla I(0),z\rangle_H|
&=\bigg|\int_{\R^n}E_1\phi+\int_{\R^n}E_2\psi\bigg|
\le C_0 \frac{\varepsilon^{(n-2)/2}}{|\ln\varepsilon|}
\big(\|\phi\|+\|\psi\|\big),
\end{aligned}
$$
ehre $E_1$ and $E_2$ were defined in \eqref{eq:errors}.
For $z\in B_\delta(0,r_0)$ and $\tau\in(0,1]$,
Lemma~\ref{lem:convexity}, applied to $\tau z$ and $0$,
gives
$$
\langle\nabla I(\tau z)-\nabla I(0),z\rangle_H
\ge\gamma\tau\|z\|_H^2.
$$
Since $I(0)=0$, integration over $\tau\in[0,1]$ yields
$$
\begin{aligned}
I(z)
&=\int_0^1\langle\nabla I(\tau z),z\rangle_H\,d\tau \ge\langle\nabla I(0),z\rangle_H
+\frac{\gamma}{2}\|z\|_H^2\\
&\ge
-C_0 \frac{\varepsilon^{(n-2)/2}}{|\ln\varepsilon|}
\big(\|\phi\|+\|\psi\|\big)
+\frac{\gamma}{4}
\big(\|\phi\|+\|\psi\|\big)^2.
\end{aligned}
$$
Choose $R>0$ sufficiently large, independently of $\beta$.
Since $\frac{\varepsilon^{(n-2)/2}}{|\ln\varepsilon|}\to0$ as $\beta\to0^-$,
we may then choose $\beta_0<0$ sufficiently close to zero
so that Lemma~\ref{lem:convexity} applies and
$$
R\left(\frac{\varepsilon^{(n-2)/2}}{|\ln\varepsilon|}
\right)<r_0
\qquad\text{for every }\beta\in[\beta_0,0).
$$
The functional $I$ is continuous and strictly convex
on the closed convex ball
$B_\delta(0,R \left(\frac{\varepsilon^{(n-2)/2}}{|\ln\varepsilon|}\right))$.
Its restriction to this ball is therefore weakly
lower semicontinuous. Since the ball is weakly compact,
$I$ attains its minimum there at a unique point
$z_*=(\phi_*,\psi_*)$.

On the boundary of the ball, the preceding lower bound gives
$$
I(z)\ge
\left(\frac{\varepsilon^{(n-2)/2}}{|\ln\varepsilon|}\right)^2
\left(\frac{\gamma R^2}{4}-C_0R\right)>0=I(0),
$$
by the choice of $R$. Thus $z_*$ lies in the interior
of the ball and satisfies $\nabla I(z_*)=0$.
This is precisely \eqref{eq:correction-eqs}.

Finally, if $z,z'\in B_\delta(0,r_0)$ are critical points
of $I$, Lemma~\ref{lem:convexity} gives
$$
\gamma\|z-z'\|_H^2
\le
\langle\nabla I(z)-\nabla I(z'),z-z'\rangle_H
=0.
$$
Hence $z=z'$, proving uniqueness among all critical
points in $B_\delta(0,r_0)$.
\end{proof}

\begin{lemma}
\label{lem:continuity}
Fix a compact set $D\subset(0,\infty)$ and $T>0$.
For $\beta<0$ sufficiently close to zero, fix
$\varepsilon$ as in \eqref{eq:eps-beta} and set
$$
\mathcal D_\beta:=
\left\{
\big(1+t\varepsilon^{n-2},
     1+s\varepsilon^{n-2},
     d\varepsilon\big):
d\in D,\ |t|+|s|\le T
\right\}.
$$
Then the map
$$
\mathcal D_\beta\ni(a,b,\delta)
\longmapsto
(\phi_{a,b,\delta},\psi_{a,b,\delta})
\in\dot H^1(\R^n)\times\dot H^1(\R^n)
$$
given by Proposition~\ref{lem:minimizer} is continuous.
Consequently, the map
defined in \eqref{eq:bifurcation-eqs} is continuous
on $\mathcal D_\beta$.
\end{lemma}

\begin{proof}
Throughout the proof, $\beta$ is fixed.
Write $\lambda=(a,b,\delta)$, $V_\delta:=V$, and
$$
z_\lambda:=
(\phi_{a,b,\delta},\psi_{a,b,\delta}).
$$
On $\mathcal X:=X_k\times X_k$, use the
Hilbert inner product and norm
$$
\langle(\phi,\psi),(h,g)\rangle_H
:=\langle\phi,h\rangle_H+\langle\psi,g\rangle_H,
\qquad
\|(\phi,\psi)\|_H:=
\big(\|\phi\|^2+\|\psi\|^2\big)^{1/2}.
$$
Set
$$
\mathcal H_\lambda:=
\mathcal H_1\times\mathcal H_2,
\qquad
\Pi_\lambda:=(P_1^\perp,P_2^\perp).
$$
Only the second factor and the second projection
depend on $\delta$.

\medskip
\noindent
\emph{Step 1: continuity of the projections.}
The bubbles $U_{\delta,\xi}$ depend smoothly on
$(\delta,\xi)\in(0,\infty)\times\R^n$
as elements of $\dot H^1(\R^n)$.
Since $\xi_j(\delta)$ is smooth for $\delta\in(0,1)$,
the maps
$$
\delta\mapsto V_\delta,
\qquad
\delta\mapsto Z_\delta
=\delta\frac{d}{d\delta}V_\delta
$$
are continuous into $\dot H^1(\R^n)$.
Since $-\Delta U_i=U_i^p$,
$$
\|V_\delta\|^2
=\sum_{i,j}\int_{\R^n}U_i^pU_j
\ge k\int_{\R^n}U^{p+1}.
$$
Let $\vec Z:=(Z^1,\ldots,Z^n)$. Directly from
\eqref{eq:Zdelta},
$$
Z_j(\xi_j+\delta y)
=\delta^{-(n-2)/2}
\left[
Z^0(y)-\frac{\delta}{\rho}\,\zeta_j\cdot\vec Z(y)
\right].
$$
By scale invariance of the norm and parity,
$$
\|Z_j\|^2
=\|Z^0\|^2
+\frac{\delta^2}{\rho^2}\|\zeta_j\cdot\vec Z\|^2
\ge\|Z^0\|^2.
$$
Moreover, differentiating the bubble equation gives
$$
-\Delta Z_i=pU_i^{p-1}Z_i.
$$
Direct differentiation of
$$
U_{\delta,\xi}(x)
=C_n\delta^{(n-2)/2}
\big(\delta^2+|x-\xi|^2\big)^{-(n-2)/2}
$$
gives
\begin{equation}\label{eq:pointwise-modes}
\begin{aligned}
\delta\,\partial_\delta U_{\delta,\xi}
&=\frac{n-2}{2}U_{\delta,\xi}
\frac{|x-\xi|^2-\delta^2}{|x-\xi|^2+\delta^2},\qquad\delta\,\nabla_\xi U_{\delta,\xi}
&=(n-2)U_{\delta,\xi}
\frac{\delta(x-\xi)}{|x-\xi|^2+\delta^2}.
\end{aligned}
\end{equation}
Here $\partial_\delta$ denotes differentiation with
$\xi$ fixed. Since
$$
\delta|x-\xi|
\le\frac12\big(|x-\xi|^2+\delta^2\big),
$$
we obtain
\begin{equation}\label{eq:modes-dominated}
\big|\delta\,\partial_\delta U_{\delta,\xi}\big|
\le\frac{n-2}{2}U_{\delta,\xi},
\qquad
\big|\delta\,\nabla_\xi U_{\delta,\xi}\big|
\le\frac{n-2}{2}U_{\delta,\xi}.
\end{equation}
Writing $\xi_j=\rho \zeta_j$, with
$\rho=\sqrt{1-\delta^2}$ and $\zeta_j$ independent of
$\delta$, the total derivative is
$$
Z_j
=\delta\,\partial_\delta U_{\delta,\xi_j}
-\frac{\delta^2}{\rho}
\zeta_j\cdot\nabla_\xi U_{\delta,\xi_j}.
$$
For $\delta$ sufficiently small, $\delta/\rho\le1$,
and therefore
\begin{equation}\label{boundZj}
|Z_j|\le(n-2)U_j.
\end{equation}
Thus, for $i\ne j$ (see \eqref{eq:eps-asympt} below) yields
$$
|\langle Z_i,Z_j\rangle_H|
\le p(n-2)^2\int_{\R^n}U_i^pU_j
=O(\delta^{n-2}),
$$
and consequently
$$
\|Z_\delta\|^2
\ge k\|Z^0\|^2-O(\delta^{n-2}).
$$
Since the norm of each bubble is independent of its
parameters,
$$
\langle U_j,Z_j\rangle_H
=\frac{\delta}{2}\frac{d}{d\delta}\|U_j\|^2=0.
$$
For $i\ne j$,
$$
|\langle U_i,Z_j\rangle_H|
=\left|\int_{\R^n}U_i^pZ_j\right|
\le(n-2)\int_{\R^n}U_i^pU_j
=O(\delta^{n-2}).
$$
Hence
$$
\langle V_\delta,Z_\delta\rangle_H
=O(\delta^{n-2}).
$$
These estimates, together with $\|V_\delta\|\le k\|U\|$,
show that
$$
\|V_\delta\|^2\|Z_\delta\|^2
-\langle V_\delta,Z_\delta\rangle_H^2\geq c_0>0,$$
for $\delta$ sufficiently small.

To express the projection, set
$$
E_1^\delta:=V_\delta,\qquad E_2^\delta:=Z_\delta,
\qquad
G_\delta:=
\big(\langle E_i^\delta,E_j^\delta\rangle_H\big)_{i,j=1}^2.
$$
Then the matrix $G_\delta$ is invertible and
$$
P_2h=
\sum_{i,j=1}^2
(G_\delta^{-1})_{ij}
\langle h,E_j^\delta\rangle_H E_i^\delta.
$$
It follows that $\lambda\mapsto\Pi_\lambda$ is continuous
in operator norm.

\medskip
\noindent
\emph{Step 2: uniform continuity of the residual map
on bounded sets.}
Let
$$
\mathcal S(u,v):=
\big((-\Delta)^{-1}R_1(u,v),
     (-\Delta)^{-1}R_2(u,v)\big).
$$
By Remark~\ref{rem:symmetric-criticality},
$\mathcal S$ maps $\mathcal X$ into $\mathcal X$.
For every $\varrho>0$, we claim that
$$
\|\mathcal S(u,v)-\mathcal S(u',v')\|_H
\le C_\varrho(\tau+\tau^q),
\qquad
\tau:=\|u-u'\|+\|v-v'\|,
$$
whenever the norms of $u,v,u',v'$ are at most $\varrho$.
By the Sobolev inequality, it suffices to estimate the nonlinear terms.
The inequalities
$$
|x_+^p-y_+^p|
\lesssim (|x|^{p-1}+|y|^{p-1})|x-y|,
$$
$$
|x_+^q-y_+^q|\le|x-y|^q,
\qquad
|x_+^{q+1}-y_+^{q+1}|
\lesssim(|x|^q+|y|^q)|x-y|
$$
hold for all $x,y\in\R$.
Since $p=2q+1$ and $p(2^*)'=2^*$,
H\"older's inequality gives
$$
\|u_+^p-(u')_+^p\|_{(2^*)'}
\lesssim
\big(\|u\|_{2^*}^{p-1}
+\|u'\|_{2^*}^{p-1}\big)
\|u-u'\|_{2^*},
$$
and
$$
\begin{aligned}
\|u_+^qv_+^{q+1}
-(u')_+^q(v')_+^{q+1}\|_{(2^*)'}
&\lesssim
\|u-u'\|_{2^*}^q\|v\|_{2^*}^{q+1}+
\|u'\|_{2^*}^q
\big(\|v\|_{2^*}^q+\|v'\|_{2^*}^q\big)
\|v-v'\|_{2^*}.
\end{aligned}
$$
The remaining terms are treated in the same way.
The Sobolev inequality proves the claimed modulus
of continuity. In particular, $\mathcal S$ is bounded
on bounded sets.

\medskip
\noindent
\emph{Step 3: continuity of the correction.}
Define
$$
\mathcal W(\lambda,z):=
\mathcal S(aU+\phi,bV_\delta+\psi),
\qquad z=(\phi,\psi)\in\mathcal X.
$$
The gradient of $I_{a,b,\delta}$ on
$\mathcal H_\lambda$ is
$\Pi_\lambda\mathcal W(\lambda,\cdot)$ (see \eqref{eq:correction-eqs}).
Thus
$$
z_\lambda\in\mathcal H_\lambda,
\qquad
\Pi_\lambda\mathcal W(\lambda,z_\lambda)=0.
$$
Let $\lambda_m\to\lambda$ in $\mathcal D_\beta$,
and write $z_m:=z_{\lambda_m}$.
Set
$$
\widehat z_m:=\Pi_\lambda z_m\in\mathcal H_\lambda.
$$
By Proposition~\ref{lem:minimizer},
$z_m\in B_{\delta_m}(0,r_0)$.
Since the projections act componentwise and do not
increase norms, also
$\widehat z_m\in B_\delta(0,r_0)$.
Moreover, by Step 1,
$$
\|\widehat z_m-z_m\|_H
=\|(\Pi_\lambda-\Pi_{\lambda_m})z_m\|_H
\le r_0\|\Pi_\lambda-\Pi_{\lambda_m}\|
\longrightarrow0.
$$
Apply Lemma~\ref{lem:convexity} at the fixed parameter
$\lambda$ to $\widehat z_m$ and $z_\lambda$.
Since $\Pi_\lambda\mathcal W(\lambda,z_\lambda)=0$ and $\Pi_{\lambda_m}\mathcal W(\lambda_m,z_m)=0$,
using the fact that $\|\Pi_\lambda\|\leq 1$ we obtain
$$
\gamma\|\widehat z_m-z_\lambda\|_H
\le
\|\Pi_\lambda\mathcal W(\lambda,\widehat z_m)\|_H\leq 
\|\mathcal W(\lambda,\widehat z_m)
-\mathcal W(\lambda_m,z_m)\|_H+
\|\Pi_\lambda-\Pi_{\lambda_m}\|
\,\|\mathcal W(\lambda_m,z_m)\|_H.
$$
The second term tends to zero by Steps 1 and 2.

For the first term, the two pairs
$$
(aU+\widehat\phi_m,bV_\delta+\widehat\psi_m),
\qquad
(a_mU+\phi_m,b_mV_{\delta_m}+\psi_m)
$$
belong to a fixed bounded subset of
$\dot H^1(\R^n)\times\dot H^1(\R^n)$.
Their distance tends to zero because
$\lambda_m\to\lambda$,
$V_{\delta_m}\to V_\delta$, and
$\widehat z_m-z_m\to0$.
The uniform continuity proved in Step 2 therefore
shows that the first term also tends to zero.

Hence $\widehat z_m\to z_\lambda$, and consequently
$z_m\to z_\lambda$ in the product Hilbert norm.
This proves the continuity of the correction map.

\medskip
\noindent
\emph{Step 4: continuity of the reduced field.}
By Steps 2 and 3, the map
$$
\lambda\longmapsto
\mathcal W(\lambda,z_\lambda)
$$
is continuous into $\mathcal X$.
The testing functions $U$,
$V$, and $Z_\delta$
also depend continuously on the parameters in
$\dot H^1(\R^n)$.
Therefore each pairing in
\eqref{eq:bifurcation-eqs} is continuous.
%
This completes the proof.
\end{proof}

\begin{remark}\label{rem:only-continuous}
Lemma~\ref{lem:continuity} proves that the correction depends continuously on the parameters.
Our argument does not establish $C^1$ dependence.
We therefore solve the finite-dimensional problem directly
through the residual equations, without requiring
differentiability of the correction map.
\end{remark}

\section{The finite-dimensional reduction}\label{sec:finite-dim}

By Proposition~\ref{lem:minimizer}, for every admissible
choice of $(a,b,\delta)$ there exists a unique small
correction $(\phi,\psi)$ satisfying
\eqref{eq:correction-eqs}.
It remains to impose the three bifurcation equations
\eqref{eq:bifurcation-eqs}.

Using the residuals defined in \eqref{eq:residuals}, set
\begin{equation}\label{eq:finite_system}
\begin{aligned}
\Lambda_1(a,b,\delta)
&:=\big\langle R_1(aU+\phi,bV+\psi),U\big\rangle,\\
\Lambda_2(a,b,\delta)
&:=\big\langle R_2(aU+\phi,bV+\psi),U_1\big\rangle,\\
\Lambda_3(a,b,\delta)
&:=\big\langle R_2(aU+\phi,bV+\psi),Z_\delta\big\rangle.
\end{aligned}
\end{equation}
We seek parameters such that
$\Lambda_1=\Lambda_2=\Lambda_3=0$. Although $U_1$ does not belong to $X_k$, testing the
second residual against $U_1$ is equivalent to testing
it against $V$. Indeed, rotational symmetry gives
\begin{equation}\label{eq:U1-vs-V}
\big\langle R_2(u,v),V\big\rangle
=k\big\langle R_2(u,v),U_1\big\rangle
\qquad\text{for all }u,v\in X_k.
\end{equation}
Thus the system $\Lambda=0$ is equivalent to
\eqref{eq:bifurcation-eqs}.

Define the constants
\begin{equation}\label{eq:interaction-constants}
\begin{gathered}
c_n:=C_n\int_{\R^n}U^p\,dx,\qquad
c_n':=C_n^{2^*}|\mathbb S^{n-1}|2^{-n/2},\\
\zeta_j:=\xi_j/\rho,\qquad
\gamma_k:=\sum_{i\ne j}|\zeta_i-\zeta_j|^{2-n},\qquad
c_k:=c_n\gamma_k,\qquad c'_k:=kc'_n.
\end{gathered}
\end{equation}
They depend only on $n$ and $k$.

Set
$$
S_0:=\int_{\R^n}U^{p+1},
$$
and define
\begin{equation}\label{eq:A-constants}
A:=\frac{n-2}{2}c_n\gamma_k,
\qquad
A':=\frac{n}{p+1}c_n'k
=\frac{n-2}{2}c_n'k,
\end{equation}
where $c_n,c_n'$ and $\gamma_k$ are defined in
\eqref{eq:interaction-constants}. These constants are positive
and depend only on $n$ and $k$.

\begin{proposition}\label{prop:Lambda-tilde}
Fix a compact set $D\subset(0,\infty)$ and $T>0$.
Let $\beta\to0^-$, let $\varepsilon$ be defined by
\eqref{eq:eps-beta}, and choose the parameters as in
\eqref{eq:rescaling}, with
$$
d\in D,\qquad |t|+|s|\le T.
$$
Then the reduced equations admit the expansions
\begin{align}
\Lambda_1(a,b,\delta)
&=\varepsilon^{n-2}
\left[-(p-1)S_0t+c_n'k\,d^{n/2}\right]
+o(\varepsilon^{n-2}),
\label{eq:L1-tilde}\\
\Lambda_2(a,b,\delta)
&=\varepsilon^{n-2}
\left[-(p-1)S_0s
-pc_n\frac{\gamma_k}{k}d^{n-2}
+c_n'd^{n/2}\right]
+o(\varepsilon^{n-2}),
\label{eq:L2-tilde}\\
\Lambda_3(a,b,\delta)
&=\varepsilon^{n-2}
\left[-Ad^{n-2}+A'd^{n/2}\right]
+o(\varepsilon^{n-2}),
\label{eq:L3-tilde}
\end{align}
uniformly on the specified parameter set.

Let $F=(F_1,F_2,F_3)$ denote the vector field whose
components are the expressions in square brackets in
\eqref{eq:L1-tilde}--\eqref{eq:L3-tilde}.
Equivalently,
\begin{equation}\label{eq:limit-map}
\Lambda\big(
1+t\varepsilon^{n-2},
1+s\varepsilon^{n-2},
d\varepsilon
\big)
=\varepsilon^{n-2}\big(F(t,s,d)+o(1)\big),
\end{equation}
uniformly on the same parameter set.
\end{proposition}

The proof of Proposition~\ref{prop:Lambda-tilde}
is given at the end of this section. A series of technical lemmas are needed first.
Let $\varepsilon$  be defined by
\eqref{eq:eps-beta}, and
$$
a=1+t\varepsilon^{n-2},
\qquad
b=1+s\varepsilon^{n-2},
\qquad
\delta=d\varepsilon,
$$
where $d$ ranges in a fixed compact set
$D\subset(0,\infty)$ and $|t|+|s|\le T$.
Let the pair $(\phi,\psi)$ denote the correction given by
Proposition~\ref{lem:minimizer}.
All estimates are uniform on this parameter set.

\begin{lemma}
\label{lem:defect-derivative}
Let
$$
V_\delta:=\sum_{j=1}^kU_{\delta,\xi_j(\delta)},
\qquad
Z_\delta:=\delta\frac{d}{d\delta}V_\delta,
$$
where the centers are given by \eqref{eq:vertices}.
If $\delta=d\varepsilon$, with $d$ ranging in a fixed
compact subset $D$ of $(0,\infty)$, then
$$
\left\|
\delta\frac{d}{d\delta}
\big(\Delta V_\delta+V_\delta^p\big)
\right\|_{(2^*)'}
\lesssim\eta_n(\varepsilon),
$$
where
$$
\eta_n(\varepsilon):=
\begin{cases}
\varepsilon^3, & n=5,\\
\varepsilon^4|\ln\varepsilon|^{2/3}, & n=6,\\
\varepsilon^{(n+2)/2}, & n\ge7.
\end{cases}
$$
The estimate is uniform for $d\in D$.
Moreover,
\begin{equation}\label{eq:approx-kernel}
-\Delta Z_\delta-pV_\delta^{p-1}Z_\delta
=-\delta\frac{d}{d\delta}
\big(\Delta V_\delta+V_\delta^p\big).
\end{equation}
\end{lemma}

\begin{proof}
For brevity, write $V=V_\delta$,
$U_j=U_{\delta,\xi_j(\delta)}$, and
$$
Z_j:=\delta\frac{d}{d\delta}
U_{\delta,\xi_j(\delta)},
\qquad
Z_\delta=\sum_{j=1}^kZ_j.
$$
By \eqref{boundZj}, $
|Z_j|\le(n-2)U_j,
$ and summing over $j$ gives
\begin{equation}\label{eq:Z-dominated}
|Z_\delta|\le(n-2)V.
\end{equation}

\medskip
\noindent
\emph{Step 1: differentiation of the defect.}
Since $-\Delta U_j=U_j^p$, set
$$
N_\delta:=\Delta V+V^p
=V^p-\sum_{j=1}^kU_j^p.
$$
Differentiating along the family gives
$$
\begin{aligned}
\delta\frac{d}{d\delta}N_\delta
&=pV^{p-1}Z_\delta
-p\sum_{j=1}^kU_j^{p-1}Z_j=p\sum_{j=1}^k
\big(V^{p-1}-U_j^{p-1}\big)Z_j.
\end{aligned}
$$
It also gives
$$
\delta\frac{d}{d\delta}N_\delta
=\Delta Z_\delta+pV^{p-1}Z_\delta,
$$
which proves \eqref{eq:approx-kernel}.

\medskip
\noindent
\emph{Step 2: pointwise estimate on $\Omega_1$.}
On the region $\Omega_1$ defined in \eqref{eq:Omega},
write
$$
W:=\sum_{j=2}^kU_j.
$$
Then
$$
W\le(k-1)U_1,
\qquad
U_1\le V\le kU_1.
$$
The mean value theorem gives
$$
V^{p-1}-U_1^{p-1}
=(p-1)\zeta^{p-2}W
$$
for some $\zeta\in[U_1,kU_1]$.
If $p\ge2$, use $\zeta\le kU_1$;
if $p<2$, use $\zeta\ge U_1$.
In either case,
$$
|V^{p-1}-U_1^{p-1}|
\lesssim U_1^{p-2}W.
$$
Thus
$$
|V^{p-1}-U_1^{p-1}|\,|Z_1|
\lesssim U_1^{p-1}W.
$$
For $j\ne1$, using $U_j\le U_1$ and $V\le kU_1$,
we have
$$
|V^{p-1}-U_j^{p-1}|\,|Z_j|
\lesssim U_1^{p-1}U_j.
$$
Summing these estimates yields
$$
\left|\delta\frac{d}{d\delta}N_\delta\right|
\lesssim U_1^{p-1}W
\qquad\text{on }\Omega_1.
$$

\medskip
\noindent
\emph{Step 3: the norm estimate.}
The derivative of the defect is invariant under
the rotations defining $X_k$. Consequently,
$$
\left\|\delta\frac{d}{d\delta}N_\delta\right\|_{(2^*)'}
\lesssim
\|U_1^{p-1}W\|_{L^{(2^*)'}(\Omega_1)}.
$$
The estimate established in the proof of
Proposition~\ref{prop:error} gives
$$
\|U_1^{p-1}W\|_{L^{(2^*)'}(\Omega_1)}
\lesssim
\begin{cases}
\delta^3, & n=5,\\
\delta^4|\ln\delta|^{2/3}, & n=6,\\
\delta^{(n+2)/2}, & n\ge7.
\end{cases}
$$
Finally, since $\delta=d\varepsilon$ and $d\in D$,
we have
$
\delta\asymp\varepsilon,
$
$
|\ln\delta|\asymp|\ln\varepsilon|
$
uniformly as $\varepsilon\to0$.
This proves the asserted estimate.
\end{proof}

\begin{lemma}
\label{lem:coupling-integral}
Let $n\ge5$ and let $k\ge2$ be fixed. 
Let $\varepsilon$ be defined by \eqref{eq:eps-beta}
and let $\delta=d\varepsilon$, with $d$ ranging
in a fixed compact subset $D$ of $(0,\infty)$.
Then
\begin{align}
|\beta|\int_{\R^n}U^{q+1}V^{q+1}\,dx
&=c'_k d^{n/2}\varepsilon^{n-2}
+O\!\left(
\frac{\varepsilon^{n-2}}{|\ln\varepsilon|}
\right),
\label{eq:coupling-integral-V}\\
|\beta|\int_{\R^n}U^{q+1}V^qU_1\,dx
&=\frac{c'_k}{k}d^{n/2}\varepsilon^{n-2}
+O\!\left(
\frac{\varepsilon^{n-2}}{|\ln\varepsilon|}
\right),
\label{eq:coupling-integral-U1}\\
|\beta|\int_{\R^n}U^{q+1}V^qZ_\delta\,dx
&=\frac{n-2}{2}c'_k d^{n/2}\varepsilon^{n-2}
+O\!\left(
\frac{\varepsilon^{n-2}}{|\ln\varepsilon|}
\right),
\label{eq:Z-coupling}
\end{align}
where $c_k'$ is defined in \eqref{eq:interaction-constants}.
All estimates are uniform for $d\in D$.
\end{lemma}

\begin{proof}
We first work with $\delta$ and then substitute
$\delta=d\varepsilon$.
Write
$$
\alpha:=q+1=\frac{n}{n-2},
\qquad
r:=|x-\xi_1|,
\qquad
W:=\sum_{j=2}^kU_j.
$$

\medskip
\noindent
\emph{Step 1: reduction of the first integral to one bubble.}
On $\Omega_1$, defined in \eqref{eq:Omega}, we have
$V=U_1+W$ and $U_j\le U_1$ for every $j$.
Since $k$ is fixed and $\rho\to1$, the mutual distances
between the centers are bounded away from zero.
Moreover, for $x\in\Omega_1$ and $j\ne1$,
$$
|x-\xi_j|\ge r,
\qquad
|x-\xi_j|\ge\frac12|\xi_1-\xi_j|.
$$
Consequently,
$$
W\lesssim\delta^{(n-2)/2}(1+r)^{2-n}.
$$
We also have
$$
U(x)^\alpha\lesssim(1+r)^{-n},
\qquad
U_1(x)\lesssim\delta^{(n-2)/2}r^{2-n},
\qquad
U_1(x)^q\lesssim\delta r^{-2}.
$$
Since $0<q<1$,
$$
0\le V^\alpha-U_1^\alpha
\le\alpha V^qW
\lesssim U_1^qW+W^\alpha.
$$
Hence
$$
U^\alpha(V^\alpha-U_1^\alpha)
\lesssim
\delta^{n/2}
\left[
r^{-2}(1+r)^{2-2n}+(1+r)^{-2n}
\right].
$$
Both functions on the right are integrable in $\R^n$.
By rotational symmetry,
$$
\int_{\R^n}U^\alpha V^\alpha\,dx
=k\int_{\Omega_1}U^\alpha U_1^\alpha\,dx
+O(\delta^{n/2}).
$$

There exists a fixed $r_*>0$ such that
$B(\xi_1,r_*)\subset\Omega_1$ for all sufficiently
small $\delta$. On $\R^n\setminus\Omega_1$,
we therefore have $r\ge r_*$ and
$$
U^\alpha U_1^\alpha
\lesssim
\delta^{n/2}(1+r)^{-n}r^{-n}.
$$
Thus
$$
\int_{\R^n\setminus\Omega_1}
U^\alpha U_1^\alpha\,dx=O(\delta^{n/2}),
$$
and
$$
\int_{\R^n}U^\alpha V^\alpha\,dx
=k\int_{\R^n}U^\alpha U_1^\alpha\,dx
+O(\delta^{n/2}).
$$

\medskip
\noindent
\emph{Step 2: the one-bubble integral.}
Changing variables $x=\xi_1+\delta y$, we obtain
$$
\int_{\R^n}U^\alpha U_1^\alpha\,dx
=\delta^{n/2}\int_{\R^n}
U(\xi_1+\delta y)^\alpha U(y)^\alpha\,dy.
$$
The region $|y|<1$ contributes $O(1)$ to the last
integral. On $1\le|y|\le\delta^{-1}$,
$$
U(y)^\alpha
=C_n^\alpha|y|^{-n}+O(|y|^{-n-2}),
$$
and
$$
U(\xi_1+\delta y)^\alpha
=U(\xi_1)^\alpha+O(\delta|y|).
$$
Since $|\xi_1|^2=1-\delta^2$,
$$
U(\xi_1)^\alpha
=C_n^\alpha(2-\delta^2)^{-n/2}
=C_n^\alpha2^{-n/2}+O(\delta^2).
$$
It follows that
$$
\begin{aligned}
&\int_{\{1\le|y|\le\delta^{-1}\}}
U(\xi_1+\delta y)^\alpha U(y)^\alpha\,dy=
C_n^{2\alpha}2^{-n/2}|\mathbb S^{n-1}|
|\ln\delta|+O(1).
\end{aligned}
$$
Indeed, the errors are bounded by a constant multiple of
$$
\int_1^{\delta^{-1}}
\left(r^{-3}+\delta+\frac{\delta^2}{r}\right)\,dr
=O(1).
$$

On $|y|>\delta^{-1}$, using
$U(y)^\alpha\lesssim|y|^{-n}$ and setting $z=\delta y$,
we get
$$
\begin{aligned}
\int_{\{|y|>\delta^{-1}\}}
U(\xi_1+\delta y)^\alpha U(y)^\alpha\,dy\lesssim
\int_{\{|z|>1\}}|z|^{-n}U(\xi_1+z)^\alpha\,dz
=O(1).
\end{aligned}
$$
This bound is uniform for $|\xi_1|\le1$:
on $1<|z|\le3$ the integrand is bounded,
and for $|z|>3$ it is $O(|z|^{-2n})$.

Since $2\alpha=2^*$, Steps 1 and 2 give
$$
\int_{\R^n}U^\alpha V^\alpha\,dx
=c'_k\delta^{n/2}|\ln\delta|
+O(\delta^{n/2}).
$$

\medskip
\noindent
\emph{Step 3: the integral tested against $U_1$.}
Rotational symmetry implies that
$$
\int_{\R^n}U^\alpha V^qU_j\,dx
$$
is independent of $j$. Summing over $j$ gives the
exact identity
$$
\int_{\R^n}U^\alpha V^qU_1\,dx
=\frac1k\int_{\R^n}U^\alpha V^\alpha\,dx.
$$

\medskip
\noindent
\emph{Step 4: the integral tested against $Z_\delta$.}
Write $Z_\delta=\sum_jZ_j$ and recall \eqref{boundZj}.
On $\Omega_1$, since $0\le V^q-U_1^q\le W^q$, we have
$$
\begin{aligned}
|V^qZ_\delta-U_1^qZ_1|
&\le(V^q-U_1^q)|Z_\delta|
+U_1^q\left|\sum_{j=2}^kZ_j\right|\lesssim W^qU_1+W^{q+1}+U_1^qW.
\end{aligned}
$$
The estimates in Step 1 yield
$$
\begin{aligned}
U^\alpha W^qU_1
&\lesssim
\delta^{n/2}r^{2-n}(1+r)^{-n-2},\\
U^\alpha W^{q+1}
&\lesssim
\delta^{n/2}(1+r)^{-2n},\\
U^\alpha U_1^qW
&\lesssim
\delta^{n/2}r^{-2}(1+r)^{2-2n}.
\end{aligned}
$$
All three functions of $r$ are integrable.
Therefore, by rotational symmetry,
$$
\int_{\R^n}U^\alpha V^qZ_\delta\,dx
=k\int_{\Omega_1}U^\alpha U_1^qZ_1\,dx
+O(\delta^{n/2}).
$$
Set $\vec Z:=(Z^1,\ldots,Z^n)$ and
$\widetilde\Omega_1:=(\Omega_1-\xi_1)/\delta$.
By the definition of $Z_1$,
$$
Z_1(\xi_1+\delta y)
=\delta^{-(n-2)/2}
\left[
Z^0(y)-\frac{\delta}{\rho}e_1\cdot\vec Z(y)
\right].
$$
Moreover,
$$
Z^0(y)
=\frac{n-2}{2}U(y)
-(n-2)\frac{U(y)}{1+|y|^2}.
$$
Changing variables $x=\xi_1+\delta y$, we obtain
$$
\begin{aligned}
\int_{\Omega_1}U^\alpha U_1^qZ_1\,dx
={}&\delta^{n/2}
\int_{\widetilde\Omega_1}
U(\xi_1+\delta y)^\alpha U(y)^q\cdot
\left[
\frac{n-2}{2}U(y)
-(n-2)\frac{U(y)}{1+|y|^2}
-\frac{\delta}{\rho}e_1\cdot\vec Z(y)
\right]\,dy.
\end{aligned}
$$
The first term equals
$$
\frac{n-2}{2}\int_{\Omega_1}
U^\alpha U_1^\alpha\,dx.
$$
The remaining terms are $O(\delta^{n/2})$ because
$U(\xi_1+\delta y)$ is uniformly bounded and
$$
\frac{U(y)^{q+1}}{1+|y|^2}
\lesssim(1+|y|)^{-n-2},
\qquad
U(y)^q|\vec Z(y)|
\lesssim(1+|y|)^{-n-1}.
$$
Hence
$$
\int_{\Omega_1}U^\alpha U_1^qZ_1\,dx
=\frac{n-2}{2}
\int_{\Omega_1}U^\alpha U_1^\alpha\,dx
+O(\delta^{n/2}).
$$
Combining this identity with Step 1 gives
$$
\int_{\R^n}U^\alpha V^qZ_\delta\,dx
=\frac{n-2}{2}
\int_{\R^n}U^\alpha V^\alpha\,dx
+O(\delta^{n/2}).
$$
Therefore,
$$
\int_{\R^n}U^\alpha V^qZ_\delta\,dx
=\frac{n-2}{2}c'_k
\delta^{n/2}|\ln\delta|
+O(\delta^{n/2}).
$$

\medskip
\noindent
\emph{Step 5: rescaling.}
For $\delta=d\varepsilon$ and $d\in D$,
$$
|\ln(d\varepsilon)|
=|\ln\varepsilon|-\ln d
=|\ln\varepsilon|+O(1)
$$
uniformly as $\varepsilon\to0$.
Substitution and multiplying by
$$
|\beta|
=\frac{\varepsilon^{(n-4)/2}}{|\ln\varepsilon|}
$$
prove \eqref{eq:coupling-integral-V}--\eqref{eq:Z-coupling}.
\end{proof}

\begin{lemma}
\label{lem:defect-integrals}
Let $\delta=d\varepsilon$, with $d$ ranging in a fixed
compact subset $D$ of $(0,\infty)$.
Then
\begin{equation}\label{eq:D1-asympt}
\int_{\R^n}
\bigg(V^p-\sum_{j=1}^kU_j^p\bigg)U_1\,dx
=\frac pk c_k d^{n-2}\varepsilon^{n-2}
+O(\varepsilon^{n-1}),
\end{equation}
and
\begin{equation}\label{eq:Z-defect}
\int_{\R^n}
\bigg(V^p-\sum_{j=1}^kU_j^p\bigg)Z_\delta\,dx
=\frac{n-2}{2}c_k d^{n-2}\varepsilon^{n-2}
+O(\varepsilon^{n-1}),
\end{equation}
where $c_k$ is defined in \eqref{eq:interaction-constants}. Both estimates are uniform for $d\in D$.
\end{lemma}

\begin{proof}
We first work with $\delta$. Write
$$
N:=V^p-\sum_{j=1}^kU_j^p,
\qquad
W:=\sum_{j=2}^kU_j,
\qquad
r:=|x-\xi_1|.
$$
Since $k$ is fixed and $\rho\to1$, we may choose
$r_*>0$, independently of $\delta$, such that
$B(\xi_1,r_*)\subset\Omega_1$ and this ball stays
a positive distance from every other center.
On $\Omega_1$,
$$
W\le(k-1)U_1,
\qquad
W\lesssim\delta^{(n-2)/2}(1+r)^{2-n},
$$
and on $B(\xi_1,r_*)$,
$$
|\nabla W|\lesssim\delta^{(n-2)/2}.
$$

\medskip
\noindent
\emph{Step 1: pointwise expansion and error estimates.}
Since $U_1\le V\le kU_1$ on $\Omega_1$,
Taylor's formula gives
$$
N=pU_1^{p-1}W
+O\big(U_1^{p-2}W^2+W^p\big).
$$
This estimate holds also when $p<2$, since the
intermediate points in Taylor's formula lie
between $U_1$ and $kU_1$.

We will use
$$
\int_{\Omega_1}U_1^{p-1}W^2\,dx=O(\delta^n),
\qquad
\int_{\Omega_1}U_1W^p\,dx=O(\delta^n).
$$
Indeed,
$$
\begin{aligned}
\int_{\Omega_1}U_1^{p-1}W^2\,dx
&\lesssim
\delta^n\int_0^\infty
\frac{r^{n-1}(1+r)^{4-2n}}
{(\delta^2+r^2)^2}\,dr \lesssim
\delta^n\int_0^\infty
r^{n-5}(1+r)^{4-2n}\,dr
\lesssim\delta^n,
\end{aligned}
$$
where $n\ge5$ ensures integrability at the origin.
Similarly,
$$
\int_{\Omega_1}U_1W^p\,dx
\lesssim
\delta^n\int_0^\infty r(1+r)^{-n-2}\,dr
\lesssim\delta^n.
$$

\medskip
\noindent
\emph{Step 2: a localization estimate.}
For every measurable function $g$ with
$|g|\le C_gU^p$, we claim that
$$
\begin{aligned}
&\int_{\Omega_1}
\delta^{-(n+2)/2}
g\left(\frac{x-\xi_1}{\delta}\right)W(x)\,dx=
\delta^{(n-2)/2}W(\xi_1)
\int_{\R^n}g(y)\,dy
+O(\delta^{n-1}),
\end{aligned}
$$
with a constant uniform when $C_g$ is bounded.

Outside $B(\xi_1,r_*)$, the absolute value of the
integrand is bounded by $C_gU_1^pW$, and
$$
\int_{\Omega_1\setminus B(\xi_1,r_*)}
U_1^pW\,dx
\lesssim
\delta^n\int_{r_*}^\infty r^{-3}\,dr
=O(\delta^n).
$$
Inside the ball, set $x=\xi_1+\delta y$.
The gradient bound for $W$ gives
$$
|W(\xi_1+\delta y)-W(\xi_1)|
\lesssim\delta^{n/2}|y|.
$$
Since $\int_{\R^n}|y|U(y)^p\,dy<\infty$,
replacing $W(\xi_1+\delta y)$ by $W(\xi_1)$
produces an error $O(\delta^{n-1})$.
Finally,
$$
\int_{\{|y|>r_*/\delta\}}U(y)^p\,dy=O(\delta^2),
$$
and $W(\xi_1)=O(\delta^{(n-2)/2})$, so extending
the resulting integral to $\R^n$ costs $O(\delta^n)$.
This proves the claim.

Moreover, since $\xi_j=\rho\zeta_j$ and
$\rho^2=1-\delta^2$,
$$
\begin{aligned}
W(\xi_1)
&=C_n\delta^{(n-2)/2}
\sum_{j=2}^k
\left(\delta^2+\rho^2|\zeta_1-\zeta_j|^2
\right)^{-(n-2)/2}=C_n\delta^{(n-2)/2}
\left(\frac{\gamma_k}{k}+O(\delta^2)\right).
\end{aligned}
$$
Here we used the rotational symmetry of the polygon.

\medskip
\noindent
\emph{Step 3: the integral tested against $U_1$.}
By Step 1,
$$
\int_{\Omega_1}NU_1\,dx
=p\int_{\Omega_1}U_1^pW\,dx+O(\delta^n).
$$
Applying Step 2 with $g=U^p$, we obtain
$$
\begin{aligned}
\int_{\Omega_1}U_1^pW\,dx
&=\delta^{(n-2)/2}W(\xi_1)
\int_{\R^n}U^p\,dx+O(\delta^{n-1})=\frac{c_k}{k}\delta^{n-2}
+O(\delta^{n-1}).
\end{aligned}
$$
It remains to estimate the integral outside $\Omega_1$.
For $j\ne1$, set $W_j:=\sum_{i\ne j}U_i$.
On $\Omega_j$,
$$
0\le N
\le (U_j+W_j)^p-U_j^p
\lesssim U_j^{p-1}W_j,
\qquad U_1\le W_j.
$$
Therefore, by rotation and Step 1,
$$
\int_{\Omega_j}NU_1\,dx
\lesssim\int_{\Omega_j}U_j^{p-1}W_j^2\,dx
=O(\delta^n).
$$
Summing over the regions gives
$$
\int_{\R^n}NU_1\,dx
=\frac pk c_k\delta^{n-2}
+O(\delta^{n-1}).
$$

\medskip
\noindent
\emph{Step 4: the integral tested against $Z_\delta$.}
The functions $N$ and $Z_\delta$ are rotationally
invariant, so
$$
\int_{\R^n}NZ_\delta\,dx
=k\int_{\Omega_1}NZ_\delta\,dx.
$$
Since
$$
|Z_1|\lesssim U_1,
\qquad
\bigg|\sum_{j=2}^kZ_j\bigg|\lesssim W,
$$
the expansion in Step 1 yields
$$
NZ_\delta
=pU_1^{p-1}Z_1W
+O\big(U_1^{p-1}W^2+U_1W^p\big)
\qquad\text{on }\Omega_1.
$$
Here we also used $W\le(k-1)U_1$ to absorb
the remaining products into the displayed error.
Consequently,
$$
\int_{\R^n}NZ_\delta\,dx
=k\int_{\Omega_1}pU_1^{p-1}Z_1W\,dx
+O(\delta^n).
$$
Let $\vec Z:=(Z^1,\ldots,Z^n)$.
The rescaled expression for $Z_1$ gives
$$
pU_1^{p-1}Z_1(x)
=\delta^{-(n+2)/2}
g_\delta\left(\frac{x-\xi_1}{\delta}\right),
$$
where
$$
g_\delta(y):=
pU(y)^{p-1}
\left[Z^0(y)-\frac{\delta}{\rho}
e_1\cdot\vec Z(y)\right].
$$
The bound $|g_\delta|\lesssim U^p$ is uniform for
$\delta$ small.
By oddness,
$$
\int_{\R^n}U^{p-1}e_1\cdot\vec Z\,dy=0.
$$
Furthermore, differentiating
$$
\int_{\R^n}U_{\delta,0}^p\,dx
=\delta^{(n-2)/2}\int_{\R^n}U^p\,dx
$$
at $\delta=1$ gives
$$
\int_{\R^n}pU^{p-1}Z^0\,dy
=\frac{n-2}{2}\int_{\R^n}U^p\,dy,
$$
and hence
$$
\int_{\R^n}g_\delta\,dy
=\frac{n-2}{2}\int_{\R^n}U^p\,dy.
$$
Applying the localization estimate from Step 2,
$$
\begin{aligned}
\int_{\Omega_1}pU_1^{p-1}Z_1W\,dx
&=\frac{n-2}{2}
\delta^{(n-2)/2}W(\xi_1)
\int_{\R^n}U^p\,dy+O(\delta^{n-1})=\frac{n-2}{2}\frac{c_k}{k}
\delta^{n-2}+O(\delta^{n-1}).
\end{aligned}
$$
We conclude that
$$
\int_{\R^n}NZ_\delta\,dx
=\frac{n-2}{2}c_k\delta^{n-2}
+O(\delta^{n-1}).
$$
Finally, substituting $\delta=d\varepsilon$
proves both assertions, uniformly for $d\in D$.
\end{proof}

\begin{lemma}
\label{lem:auxiliary-interactions}
Let $n\ge5$ and let $k\ge2$ be fixed.
Let $\delta=d\varepsilon$, with $d$ ranging in a fixed
compact subset $D$ of $(0,\infty)$.
Then, for $i\ne j$,
\begin{equation}\label{eq:eps-asympt}
\int_{\R^n}U_i^pU_j\,dx=O(\varepsilon^{n-2}).
\end{equation}
Moreover,
\begin{equation}\label{eq:Z-Vp}
\int_{\R^n}V^pZ_\delta\,dx
=O(\varepsilon^{n-2}).
\end{equation}
All estimates are uniform for $d\in D$.
\end{lemma}

\begin{proof}
We first prove the estimates in terms of $\delta$.
Since $k$ is fixed and $\rho\to1$, the mutual distances
between the centers are bounded away from zero.
Choose $r_*>0$, independently of $\delta$, such that
$B(\xi_i,r_*)$ stays a positive distance from $\xi_j$
whenever $i\ne j$.

On $B(\xi_i,r_*)$, we have
$U_j\lesssim\delta^{(n-2)/2}$, and therefore
$$
\begin{aligned}
\int_{B(\xi_i,r_*)}U_i^pU_j\,dx
&\lesssim
\delta^{(n-2)/2}\int_{\R^n}U_i^p\,dx=
\delta^{n-2}\int_{\R^n}U^p\,dx
=O(\delta^{n-2}).
\end{aligned}
$$
Outside $B(\xi_i,r_*)$,
$$
U_i^p
\lesssim
\delta^{(n+2)/2}(1+|x-\xi_i|)^{-n-2},
\qquad
U_j\lesssim
\delta^{(n-2)/2}|x-\xi_j|^{2-n}.
$$
Hence
$$
\begin{aligned}
\int_{\R^n\setminus B(\xi_i,r_*)}U_i^pU_j\,dx
&\lesssim
\delta^n\int_{\R^n}
(1+|x-\xi_i|)^{-n-2}|x-\xi_j|^{2-n}\,dx=O(\delta^n).
\end{aligned}
$$
The last integral is uniformly bounded:
the singularity at $\xi_j$ is integrable,
and the integrand decays as $|x|^{-2n}$ at infinity.
This proves
$$
\int_{\R^n}U_i^pU_j\,dx=O(\delta^{n-2}).
$$

To prove the second assertion, set
$$
N:=V^p-\sum_{j=1}^kU_j^p.
$$
Since $-\Delta V=\sum_jU_j^p$, then
$$
\int_{\R^n}V^pZ_\delta\,dx
=\langle V,Z_\delta\rangle_H
+\int_{\R^n}NZ_\delta\,dx.
$$
The second term is $O(\delta^{n-2})$ by
Lemma~\ref{lem:defect-integrals}.
For the first term, write
$$
\langle V,Z_\delta\rangle_H
=\sum_{i,j=1}^k\langle U_i,Z_j\rangle_H.
$$
The diagonal terms vanish because the norm of a bubble
does not depend on its parameters:
$$
\langle U_j,Z_j\rangle_H
=\frac{\delta}{2}\frac{d}{d\delta}\|U_j\|^2=0.
$$
For $i\ne j$, the bound $|Z_j|\le(n-2)U_j$ gives
$$
|\langle U_i,Z_j\rangle_H|
=\left|\int_{\R^n}U_i^pZ_j\,dx\right|
\le(n-2)\int_{\R^n}U_i^pU_j\,dx
=O(\delta^{n-2}).
$$
Consequently,
$
\langle V,Z_\delta\rangle_H=O(\delta^{n-2}),
$
and hence
$$
\int_{\R^n}V^pZ_\delta\,dx=O(\delta^{n-2}).
$$
Substituting $\delta=d\varepsilon$ completes the proof.
\end{proof}

\begin{lemma}
\label{lem:coupling}
Fix a compact set $D\subset(0,\infty)$ and $T>0$.
Let the parameters be chosen as in \eqref{eq:rescaling},
with $\varepsilon$ defined by \eqref{eq:eps-beta},
$d\in D$ and $|t|+|s|\le T$.
Let $(\phi,\psi)$ be the correction given by
Proposition~\ref{lem:minimizer}, and let $C_h>0$.

\begin{enumerate}
\item[\textup{(i)}]
If $h$ is measurable and $|h|\le C_hV$, then
$$
\begin{aligned}
\left|
\int_{\R^n}
\Big[
(aU+\phi)_+^{q+1}(bV+\psi)_+^q
-(aU)^{q+1}(bV)^q
\Big]h\,dx
\right|
\lesssim
\varepsilon(\|\phi\|+\|\psi\|).
\end{aligned}
$$

\item[\textup{(ii)}]
If $h$ is measurable and $|h|\le C_hU$, then
$$
\begin{aligned}
\left|
\int_{\R^n}
\Big[
(aU+\phi)_+^q(bV+\psi)_+^{q+1}
-(aU)^q(bV)^{q+1}
\Big]h\,dx
\right|
\lesssim
\varepsilon(\|\phi\|+\|\psi\|).
\end{aligned}
$$
\end{enumerate}
The implicit constants may depend on $C_h$ and the
fixed parameter set, but are independent of $\beta$.
\end{lemma}

\begin{proof}
For $\beta$ sufficiently close to zero,
$a,b\in[1/2,2]$.
By Proposition~\ref{lem:minimizer},
$$
\|\phi\|+\|\psi\|
\lesssim
\frac{\varepsilon^{(n-2)/2}}{|\ln\varepsilon|}.
$$
Since $q=2/(n-2)$, it follows that
$$
(\|\phi\|+\|\psi\|)^q
\lesssim
\frac{\varepsilon}{|\ln\varepsilon|^q}
\lesssim\varepsilon.
$$
We first record two elementary inequalities.
For $A>0$ and $h\in\R$,
$$
\big|(A+h)_+^{q+1}-A^{q+1}\big|
\lesssim A^q|h|+|h|^{q+1}.
$$
Moreover,
$$
\big|(A+h)_+^q-A^q\big|
\le A^{q-1}|h|.
$$
To verify the second inequality, set
$z=(A+h)_+/A$.
For every $z\ge0$ and $0<q<1$,
$$
|z^q-1|\le|z-1|,
$$
because $z\le z^q\le1$ when $0\le z\le1$,
and $1\le z^q\le z$ when $z\ge1$.
Also, $
|(A+h)_+-A|\le|h|,
$
which proves the assertion.

\medskip
\noindent
\emph{Proof of (i).}
Set $s_0:=aU$ and $t_0:=bV$.
We write
$$
\begin{aligned}
&(s_0+\phi)_+^{q+1}(t_0+\psi)_+^q
-s_0^{q+1}t_0^q=
\big[(s_0+\phi)_+^{q+1}-s_0^{q+1}\big]
(t_0+\psi)_+^q
+s_0^{q+1}
\big[(t_0+\psi)_+^q-t_0^q\big].
\end{aligned}
$$
Using the preceding inequalities,
$(t_0+\psi)_+^q\le t_0^q+|\psi|^q$,
and $|h|\le C_hV$, 
$$
\bigg|\Big[
(aU+\phi)_+^{q+1}(bV+\psi)_+^q
-(aU)^{q+1}(bV)^q
\Big]h\bigg|\leq C\big(U^q|\phi|+|\phi|^{q+1}\big)
\big(V^q+|\psi|^q\big)V
+U^{q+1}V^q|\psi|.
$$
Here the factor $V$ from the test function compensates
the factor $V^{q-1}$ in the second difference.

By Proposition~\ref{prop:error},
$$
\int_{\R^n}U^qV^{q+1}|\phi|\,dx
\lesssim\varepsilon\|\phi\|,
\qquad
\int_{\R^n}U^{q+1}V^q|\psi|\,dx
\lesssim\varepsilon\|\psi\|.
$$
For the remaining terms, H\"older's and Sobolev's
inequalities give
$$
\int_{\R^n}U^qV|\psi|^q|\phi|\,dx
\lesssim\|\psi\|^q\|\phi\|
\lesssim\varepsilon\|\phi\|,
$$
$$
\int_{\R^n}|\phi|^{q+1}V^{q+1}\,dx
\lesssim\|\phi\|^{q+1}
\lesssim\varepsilon\|\phi\|,
$$
and
$$
\int_{\R^n}|\phi|^{q+1}|\psi|^qV\,dx
\lesssim\|\phi\|^{q+1}\|\psi\|^q
\lesssim\varepsilon^2\|\phi\|.
$$
All applications of H\"older's inequality are valid
because the total homogeneity is $2q+2=2^*$,
and $\|U\|_{2^*}$ and $\|V\|_{2^*}$ are uniformly
bounded. This proves (i).

\medskip
\noindent
\emph{Proof of (ii).}
We use the decomposition
$$
\begin{aligned}
&(s_0+\phi)_+^q(t_0+\psi)_+^{q+1}
-s_0^qt_0^{q+1}=
\big[(t_0+\psi)_+^{q+1}-t_0^{q+1}\big]
(s_0+\phi)_+^q
+t_0^{q+1}
\big[(s_0+\phi)_+^q-s_0^q\big].
\end{aligned}
$$
Since $|h|\le C_hU$, 
$$
\bigg|\Big[
(aU+\phi)_+^q(bV+\psi)_+^{q+1}
-(aU)^q(bV)^{q+1}
\Big]h\bigg|\leq C \big(V^q|\psi|+|\psi|^{q+1}\big)
\big(U^q+|\phi|^q\big)U
+V^{q+1}U^q|\phi|.
$$
The two terms containing only one correction satisfy
$$
\int_{\R^n}U^{q+1}V^q|\psi|\,dx
\lesssim\varepsilon\|\psi\|,
\qquad
\int_{\R^n}U^qV^{q+1}|\phi|\,dx
\lesssim\varepsilon\|\phi\|.
$$
The other terms are estimated by
$$
\int_{\R^n}V^qU|\phi|^q|\psi|\,dx
\lesssim\|\phi\|^q\|\psi\|
\lesssim\varepsilon\|\psi\|,
$$
$$
\int_{\R^n}|\psi|^{q+1}U^{q+1}\,dx
\lesssim\|\psi\|^{q+1}
\lesssim\varepsilon\|\psi\|,
$$
and
$$
\int_{\R^n}|\psi|^{q+1}|\phi|^qU\,dx
\lesssim\|\psi\|^{q+1}\|\phi\|^q
\lesssim\varepsilon^2\|\psi\|.
$$
This proves (ii).
\end{proof}

\begin{lemma}\label{lem:positive-part}
Let $\varepsilon$ and $(a,b,\delta)$ be given by
\eqref{eq:eps-beta} and \eqref{eq:rescaling}, with $d$ in a fixed
compact set $D\subset(0,\infty)$ and $|t|+|s|\le T$.
Let $(\phi,\psi)\in\mathcal H_1\times\mathcal H_2$ be the correction
provided by Proposition~\ref{lem:minimizer}. Set $S_0:=\int_{\R^n}U^{p+1}\,dx$. Then, as $\beta\to0^-$,
\begin{equation}\label{eqA1}
\int_{\R^n}(aU+\phi)_+^pU\,dx
=S_0+ptS_0\varepsilon^{n-2}+O\!\left(\frac{\varepsilon^{n-2}}{|\ln\varepsilon|^2}\right),
\end{equation}
\begin{equation}\label{eqA2}
\int_{\R^n}(bV+\psi)_+^pU_1\,dx
=S_0+\left[psS_0+\frac{p+1}{k}c_kd^{n-2}\right]\varepsilon^{n-2}+O\!\left(\frac{\varepsilon^{n-2}}{|\ln\varepsilon|^2}\right).
\end{equation}
All estimates are uniform on the specified parameter set.
\end{lemma}

\begin{proof}
For small $|\beta|$, $a,b\in[1/2,2]$.
We first estimate the Taylor remainder, using the choices
$$
(\Theta,c,\eta,H)=(U,a,\phi,U)
\quad\text{or}\quad
(\Theta,c,\eta,H)=(V,b,\psi,U_1).
$$
In both cases $0\le H\le\Theta$ and $\|\Theta\|_{2^*}$ is uniformly
bounded. Set
$$
\varrho:=(c\Theta+\eta)_+^p-(c\Theta)^p-p(c\Theta)^{p-1}\eta.
$$
On $\{|\eta|\le c\Theta/2\}$, Taylor's formula gives
$|\varrho|\lesssim\Theta^{p-2}|\eta|^2$.
Indeed, the intermediate point lies between $c\Theta/2$ and
$3c\Theta/2$, which also controls the negative exponent $p-2$ when
$n\ge7$. On the complementary set, $\Theta\le4|\eta|$ and
$|\varrho|\lesssim|\eta|^p$. Therefore, by H\"older and Sobolev,
\begin{equation}\label{eq:taylor}
\begin{aligned}
\left|\int_{\R^n}\varrho H\,dx\right|
&\lesssim\int_{\R^n}\Theta^{p-1}|\eta|^2\,dx
+\int_{\R^n}|\eta|^{p+1}\,dx\lesssim\|\Theta\|_{2^*}^{p-1}\|\eta\|^2+\|\eta\|^{2^*}
=O(\|\eta\|^2).
\end{aligned}
\end{equation}
Here $(p-1)n/2=2^*$ and $2^*>2$.

For the first identity, the linear term vanishes because
$$
\int_{\R^n}U^p\phi\,dx=\langle U,\phi\rangle_H=0.
$$
It follows that
$$
\int_{\R^n}(aU+\phi)_+^pU\,dx=a^pS_0+O(\|\phi\|^2).
$$
Since $a=1+t\varepsilon^{n-2}$ and $t$ is bounded,
$$
a^p=1+pt\varepsilon^{n-2}+O(\varepsilon^{2(n-2)}).
$$
Proposition~\ref{lem:minimizer} allows us to conclude \eqref{eqA1}.

For the second identity, rotational symmetry and $\psi\in\mathcal H_2$ give
$$
\int_{\R^n}U_1^p\psi\,dx
=\langle U_1,\psi\rangle_H
=\frac1k\langle V,\psi\rangle_H=0.
$$
Consequently,
$$
\int_{\R^n}V^{p-1}U_1\psi\,dx
=\int_{\R^n}(V^{p-1}-U_1^{p-1})U_1\psi\,dx.
$$
Write $W_j:=\sum_{i\ne j}U_i$, and $W:=W_1$.
On $\Omega_1$, where $U_1\le V\le kU_1$, the mean value theorem yields
$$
|V^{p-1}-U_1^{p-1}|U_1\lesssim U_1^{p-1}W.
$$
On $\Omega_j$, $j\ne1$, we have $U_1\le W_j$ and $V\le kU_j$, so
$$
|V^{p-1}-U_1^{p-1}|U_1\lesssim U_j^{p-1}W_j.
$$
By rotation and the estimate proved in Proposition~\ref{prop:error},
$$
\big\|(V^{p-1}-U_1^{p-1})U_1\big\|_{(2^*)'}
\lesssim\|U_1^{p-1}W\|_{L^{(2^*)'}(\Omega_1)}
\lesssim\eta_n(\varepsilon).
$$
H\"older and Sobolev therefore bound the linear term by
$C\eta_n(\varepsilon)\|\psi\|$.
Together with \eqref{eq:taylor}, this gives
$$
\int_{\R^n}(bV+\psi)_+^pU_1\,dx
=b^p\int_{\R^n}V^pU_1\,dx
+O\big(\eta_n(\varepsilon)\|\psi\|+\|\psi\|^2\big).
$$
Finally, for every $n\ge5$,
$$
\eta_n(\varepsilon)
=o\!\left(\frac{\varepsilon^{(n-2)/2}}{|\ln\varepsilon|}\right),
\qquad
\|\phi\|+\|\psi\|
=O\!\left(\frac{\varepsilon^{(n-2)/2}}{|\ln\varepsilon|}\right).
$$
Thus the error just estimated is
$O(\varepsilon^{n-2}/|\ln\varepsilon|^2)$.
Since
$\int U_j^pU_1=\int U_1^pU_j$, we have the exact identity
$$
\int_{\R^n}V^pU_1\,dx
=S_0+\sum_{j=2}^k\int_{\R^n}U_1^pU_j\,dx
+\int_{\R^n}\bigg(V^p-\sum_{j=1}^kU_j^p\bigg)U_1\,dx.
$$
The localization calculation in the proof of
Lemma~\ref{lem:defect-integrals} gives
$$
\int_{\Omega_1}U_1^pW\,dx
=\frac{c_k}{k}\delta^{n-2}+O(\delta^{n-1}).
$$
The tail estimate in the proof of Lemma~\ref{lem:auxiliary-interactions}
also gives 
$$\int_{\R^n\setminus\Omega_1}U_1^pW\,dx=O(\delta^n),$$
since this region stays a fixed positive distance from $\xi_1$.
Therefore
$$
\sum_{j=2}^k\int_{\R^n}U_1^pU_j\,dx
=\frac{c_k}{k}d^{n-2}\varepsilon^{n-2}+O(\varepsilon^{n-1}).
$$
Combining this with \eqref{eq:D1-asympt}, we obtain
$$
\int_{\R^n}V^pU_1\,dx
=S_0+\frac{p+1}{k}c_kd^{n-2}\varepsilon^{n-2}
+O(\varepsilon^{n-1}).
$$
Moreover, $b=1+s\varepsilon^{n-2}$ gives
$$
b^p=1+ps\varepsilon^{n-2}+O(\varepsilon^{2(n-2)}).
$$
Multiplying these expansions proves \eqref{eqA2}, since
$$
\varepsilon^{n-1}+\varepsilon^{2(n-2)}
=o\!\left(\frac{\varepsilon^{n-2}}{|\ln\varepsilon|^2}\right).
$$
All constants are uniform on the prescribed parameter set.
\end{proof}

\begin{lemma}
\label{lem:freezing}
Under the assumptions of Proposition \ref{prop:Lambda-tilde},
\begin{equation}\label{eq:L1-expansion}
\begin{split}
\Lambda_1
&=-(p-1)S_0t\,\varepsilon^{n-2}
+\frac{\varepsilon^{(n-4)/2}}{|\ln\varepsilon|}
\int_{\R^n}U^{q+1}V^{q+1}
+O\!\left(
\frac{\varepsilon^{n-2}}{|\ln\varepsilon|^2}
\right),\\
\Lambda_2
&=-(p-1)S_0s\,\varepsilon^{n-2}
-\int_{\R^n}
\bigg(V^p-\sum_{j=1}^kU_j^p\bigg)U_1+\frac{\varepsilon^{(n-4)/2}}{|\ln\varepsilon|}
\int_{\R^n}U^{q+1}V^qU_1
+O\!\left(
\frac{\varepsilon^{n-2}}{|\ln\varepsilon|^2}
\right),\\
\Lambda_3
&=-\int_{\R^n}
\bigg(V^p-\sum_{j=1}^kU_j^p\bigg)Z_\delta+\frac{\varepsilon^{(n-4)/2}}{|\ln\varepsilon|}
\int_{\R^n}U^{q+1}V^qZ_\delta
+O\!\left(
\frac{\varepsilon^{n-2}}{|\ln\varepsilon|^2}
\right).
\end{split}
\end{equation}
Here $S_0=\int_{\R^n}U^{p+1}$, and all remainders
are uniform on the specified parameter set.
\end{lemma}

\begin{proof}
By Proposition~\ref{lem:minimizer},
$$
\|\phi\|+\|\psi\|
\lesssim
\frac{\varepsilon^{(n-2)/2}}{|\ln\varepsilon|}.
$$
Moreover, \eqref{eq:eps-beta} gives
$$
|\beta|
=\frac{\varepsilon^{(n-4)/2}}{|\ln\varepsilon|},
\qquad
|\beta|\delta
=d\,\frac{\varepsilon^{(n-2)/2}}{|\ln\varepsilon|}.
$$
Consequently,
$$
\|\phi\|^2+\|\psi\|^2
+|\beta|\delta(\|\phi\|+\|\psi\|)
\lesssim
\frac{\varepsilon^{n-2}}{|\ln\varepsilon|^2},
$$
and
$$
(\|\phi\|+\|\psi\|)^q
\lesssim
\frac{\varepsilon}{|\ln\varepsilon|^q}
\lesssim\delta.
$$

\medskip
\noindent
\emph{Step 1: removal of the corrections.}
Since $\langle\phi,U\rangle_H=0$, \eqref{eqA1} and part (ii) of
Lemma~\ref{lem:coupling} give
$$
\begin{aligned}
\Lambda_1
={}&(a-a^p)S_0
-\beta a^qb^{q+1}
\int_{\R^n}U^{q+1}V^{q+1}+O\!\left(
\frac{\varepsilon^{n-2}}{|\ln\varepsilon|^2}
\right).
\end{aligned}
$$
For the second equation, rotational symmetry implies
$
\langle\psi,U_1\rangle_H
=\frac1k\langle\psi,V\rangle_H=0.
$
Also,
$$
\langle V,U_1\rangle_H
=S_0+\sum_{j=2}^k\int_{\R^n}U_1^pU_j.
$$
Then \eqref{eqA2} gives
$$
\int_{\R^n}(bV+\psi)_+^pU_1
=b^p\int_{\R^n}V^pU_1
+O\!\left(\frac{\varepsilon^{n-2}}{|\ln\varepsilon|^2}\right).
$$
Using the fact that 
$
\int_{\R^n}U_j^pU_1
=\int_{\R^n}U_1^pU_j
$
and part (i) of Lemma~\ref{lem:coupling}, with
$h=U_1\le V$, we obtain
$$
\begin{aligned}
\Lambda_2
={}&(b-b^p)
\bigg(
S_0+\sum_{j=2}^k\int_{\R^n}U_1^pU_j
\bigg)-b^p\int_{\R^n}
\left(V^p-\sum_{j=1}^kU_j^p\right)U_1-\beta a^{q+1}b^q
\int_{\R^n}U^{q+1}V^qU_1
+O\!\bigg(
\frac{\varepsilon^{n-2}}{|\ln\varepsilon|^2}
\bigg).
\end{aligned}
$$
For the third equation, use
$\langle\psi,Z_\delta\rangle_H=0$ and
$|Z_\delta|\le(n-2)V$, as proved in
\eqref{eq:Z-dominated}. We have
$$
\begin{aligned}
\Lambda_3-\langle R_2(aU,bV),Z_\delta\rangle
=&-\int_{\R^n}
\big[(bV+\psi)_+^p-(bV)^p\big]Z_\delta\\
&-\beta\int_{\R^n}
\big[
(aU+\phi)_+^{q+1}(bV+\psi)_+^q
-(aU)^{q+1}(bV)^q
\big]Z_\delta.
\end{aligned}
$$
The coupling term is bounded by
$$
C|\beta|\delta(\|\phi\|+\|\psi\|)
\lesssim
\frac{\varepsilon^{n-2}}{|\ln\varepsilon|^2}
$$
by Lemma~\ref{lem:coupling}.
Write the scalar term as
$$
(bV+\psi)_+^p-(bV)^p
=pb^{p-1}V^{p-1}\psi+\varrho,
$$
where
$$
\varrho:=
(bV+\psi)_+^p-(bV)^p-p(bV)^{p-1}\psi.
$$
On $\{|\psi|\le bV/2\}$, Taylor's formula gives
$|\varrho|\lesssim V^{p-2}\psi^2$.
On the complementary set,
$V\lesssim|\psi|$ and $|\varrho|\lesssim|\psi|^p$.
Hence
$$
\begin{aligned}
\left|\int_{\R^n}\varrho Z_\delta\right|\lesssim
\int_{\R^n}V^{p-1}\psi^2
+\int_{\R^n}|\psi|^{p+1}\lesssim\|\psi\|^2+\|\psi\|^{2^*}
\lesssim\|\psi\|^2.
\end{aligned}
$$
Differentiating the defect along the family $V_\delta$
gives
$$
pV^{p-1}Z_\delta
=-\Delta Z_\delta
+\delta\frac{d}{d\delta}(\Delta V+V^p).
$$
Therefore, by orthogonality and
Lemma~\ref{lem:defect-derivative},
$$
\begin{aligned}
\left|p\int_{\R^n}V^{p-1}Z_\delta\psi\right|
&=
\left|
\int_{\R^n}
\delta\frac{d}{d\delta}(\Delta V+V^p)\,\psi
\right| \lesssim\eta_n(\varepsilon)\|\psi\|=o\!\left(
\frac{\varepsilon^{n-2}}{|\ln\varepsilon|^2}
\right),
\end{aligned}
$$
where $\eta_n$ is defined in Proposition~\ref{prop:error}.
We conclude that
$$
\Lambda_3
=\langle R_2(aU,bV),Z_\delta\rangle
+O\!\left(
\frac{\varepsilon^{n-2}}{|\ln\varepsilon|^2}
\right).
$$
Using $-\Delta V=\sum_jU_j^p$, we also have the exact
identity
\begin{equation}\label{eq:frozen-identity}
\begin{aligned}
\langle R_2(aU,bV),Z_\delta\rangle
=-b\int_{\R^n}
\bigg(V^p-\sum_{j=1}^kU_j^p\bigg)Z_\delta+(b-b^p)\int_{\R^n}V^pZ_\delta-\beta a^{q+1}b^q
\int_{\R^n}U^{q+1}V^qZ_\delta.
\end{aligned}
\end{equation}

\medskip
\noindent
\emph{Step 2: expansion of the amplitudes.}
Uniformly for $|t|+|s|\le T$, Taylor's formula gives
$$
a-a^p
=-(p-1)t\varepsilon^{n-2}
+O(\varepsilon^{2n-4}),
$$
$$
b-b^p
=-(p-1)s\varepsilon^{n-2}
+O(\varepsilon^{2n-4}),
$$
and
$$
b=1+O(\varepsilon^{n-2}),
\qquad
b^p=1+O(\varepsilon^{n-2}),
$$
$$
a^qb^{q+1}=1+O(\varepsilon^{n-2}),
\qquad
a^{q+1}b^q=1+O(\varepsilon^{n-2}).
$$
The estimates in
Lemmas~\ref{lem:coupling-integral},
\ref{lem:defect-integrals} and
\ref{lem:auxiliary-interactions} give
$$
\sum_{j=2}^k\int_{\R^n}U_1^pU_j
=O(\varepsilon^{n-2}), \qquad
\int_{\R^n}
\bigg(V^p-\sum_{j=1}^kU_j^p\bigg)U_1
=O(\varepsilon^{n-2}),
$$
$$
\int_{\R^n}
\bigg(V^p-\sum_{j=1}^kU_j^p\bigg)Z_\delta
=O(\varepsilon^{n-2}),
\qquad
\int_{\R^n}V^pZ_\delta
=O(\varepsilon^{n-2}).
$$
They also imply
$$
\begin{aligned}
|\beta|\bigg(
&\int_{\R^n}U^{q+1}V^{q+1}
+\int_{\R^n}U^{q+1}V^qU_1+\left|\int_{\R^n}U^{q+1}V^qZ_\delta\right|
\bigg)
=O(\varepsilon^{n-2}).
\end{aligned}
$$
Putting all the estimates together and substituting
$$
-\beta=|\beta|
=\frac{\varepsilon^{(n-4)/2}}{|\ln\varepsilon|}
$$
proves \eqref{eq:L1-expansion}.
\end{proof}

\begin{proof}[Proof of Proposition~\ref{prop:Lambda-tilde}]
Combining Lemmas~\ref{lem:freezing},
\ref{lem:coupling-integral} and
\ref{lem:defect-integrals}, we obtain
$$
\begin{aligned}
\Lambda_1
&=\varepsilon^{n-2}
\left[-(p-1)S_0t+c'_k d^{n/2}\right]
+O\!\left(
\frac{\varepsilon^{n-2}}{|\ln\varepsilon|}
\right),\\
\Lambda_2
&=\varepsilon^{n-2}
\left[
-(p-1)S_0s-\frac pk c_k d^{n-2}
+\frac{c'_k}{k}d^{n/2}
\right]
+O\!\left(
\frac{\varepsilon^{n-2}}{|\ln\varepsilon|}
\right),\\
\Lambda_3
&=\frac{n-2}{2}\varepsilon^{n-2}
\left[-c_k d^{n-2}+c'_k d^{n/2}\right]
+O\!\left(
\frac{\varepsilon^{n-2}}{|\ln\varepsilon|}
\right),
\end{aligned}
$$
with uniform estimates on the prescribed parameter set.
Since
$$
c_k=c_n\gamma_k,\qquad c'_k=kc'_n,
\qquad
A=\frac{n-2}{2}c_k,\qquad
A'=\frac{n-2}{2}c'_k,
$$
these are precisely the expansions asserted in
Proposition~\ref{prop:Lambda-tilde}.
\end{proof}

\section{Proof of Theorem \ref{thm:main-m2}}
The vector field $F$ in
Proposition~\ref{prop:Lambda-tilde} has a zero
$x_*=(t_*,s_*,d_*)$, where
$$
t_*:=\frac{c_n\gamma_k}{(p-1)S_0}d_*^{n-2},
\qquad
s_*:=-\frac{c_n\gamma_k}{kS_0}d_*^{n-2}.
$$
Indeed, $d_*$ is uniquely determined by the identity
\begin{equation}\label{dvalue}
c_n'k\,d_*^{n/2}
=c_n\gamma_k\,d_*^{n-2},
\end{equation}
and substitution into $F_1,F_2,F_3$ yields $F(x_*)=0$.
The Jacobian is nondegenerate, since
$$
\det DF(x_*)
=-(p-1)^2S_0^2\frac{n-4}{2}
A\,d_*^{n-3}\ne0.
$$
Choose a compact interval $D\subset(0,\infty)$
with $d_*\in\operatorname{int}D$, and choose $T>0$
such that $|t_*|+|s_*|<T$.
Then $x_*$ lies in the interior of the admissible
parameter set. Choose a sufficiently small open ball
$B$ centered at $x_*$ whose closure is contained
in this set and such that $x_*$ is the only zero
of $F$ in $\overline B$.

For $\beta<0$ sufficiently close to zero, define
$$
F_\beta(t,s,d):=
\varepsilon^{-(n-2)}
\Lambda\big(
1+t\varepsilon^{n-2},
1+s\varepsilon^{n-2},
d\varepsilon
\big).
$$
By Lemma~\ref{lem:continuity}, $F_\beta$ is continuous.
By Proposition~\ref{prop:Lambda-tilde},
$F_\beta\to F$ uniformly on $\overline B$.
Therefore, for $\beta$ sufficiently close to zero,
$$
\deg(F_\beta,B,0)
=\deg(F,B,0)
=\operatorname{sign}\det DF(x_*)\ne0,
$$
and hence $F_\beta$ has a zero
$x_\beta=(t_\beta,s_\beta,d_\beta)\in B$.

Since $F_\beta\to F$ uniformly on $\overline B$
and $x_*$ is the only zero of $F$ there,
we have $x_\beta\to x_*$.
Set
$$
a_\beta:=1+t_\beta\varepsilon^{n-2},
\qquad
b_\beta:=1+s_\beta\varepsilon^{n-2},
\qquad
\delta_\beta:=d_\beta\varepsilon,
$$
and let $(\phi_\beta,\psi_\beta)$ be the corresponding
correction given by Proposition~\ref{lem:minimizer}.
The corrected ansatz satisfies
\eqref{eq:correction-eqs} and $\Lambda=0$.
By \eqref{eq:U1-vs-V}, the residuals therefore vanish
on the finite-dimensional spaces as well as on
their orthogonal complements in $X_k$.
Remark~\ref{rem:symmetric-criticality} then shows
that $(u_\beta,v_\beta)$ is a weak solution of
\eqref{eq:positive-system} in the full space. 
Testing the equations with the negative parts gives
$u_\beta,v_\beta\ge0$.

Finally, $x_\beta\to x_*$ yields
\begin{equation}\label{eq:amplitudes}
\begin{aligned}
a_\beta
&=1+\frac{c_n\gamma_k}{(p-1)S_0}
d_*^{n-2}\varepsilon^{n-2}
+o(\varepsilon^{n-2}),\\
b_\beta
&=1-\frac{c_n\gamma_k}{kS_0}
d_*^{n-2}\varepsilon^{n-2}
+o(\varepsilon^{n-2}),\\
\delta_\beta
&=d_*\varepsilon+o(\varepsilon),
\end{aligned}
\end{equation} 
while the estimate of the
corrections follows from
Proposition~\ref{lem:minimizer}. This completes the proof of Theorem \ref{thm:main-m2}.

\section*{Acknowledgments}
I. Ceresa Dussel is financially supported by MICIU/AEI/10.13039/501100011033 (Spain), through the grant CNS2024-154515. M. Medina is supported by the grant PID2023-149451NA-I00 of MCIN/AEI/10.13039/ 501100011033/FEDER, UE. Angela Pistoia is supported by INDAM-GNAMPA project.

\section*{AI declaration}
The authors used artificial intelligence tools as auxiliary aids for assistance
with some calculations and with the drafting and language editing of
selected parts of the manuscript. All original research ideas and
scientific contributions are entirely those of the authors, who take
full responsibility for the mathematical arguments and the content
of this work.

\bibliographystyle{amsplain-initials}
\bibliography{refs}

\end{document}